\documentclass[12pt,a4paper]{amsart}
\usepackage{amssymb,amsfonts,amsthm,amsmath,graphicx,xcolor,csquotes}
\usepackage[scr]{rsfso}
\usepackage{chngcntr}
\usepackage{enumitem}

\theoremstyle{plain}
\newtheorem{theorem}{Theorem}
\newtheorem{lemma}{Lemma}

\newtheorem{corollary}{Corollary}
\newtheorem{statement}{Statement}

\numberwithin{theorem}{section}
\numberwithin{equation}{section}
\numberwithin{statement}{section}
\numberwithin{lemma}{section}
\numberwithin{definition}{section}
\numberwithin{conjecture}{section}
\numberwithin{corollary}{section}
\begin{document}
\title[Visible Point Vector Identities]{Visible Point Vector Partition Identities for Hyperpyramid Lattices}
\author{Geoffrey B Campbell}
\address{Mathematical Sciences Institute,
         The Australian National University,
         Canberra, ACT, 0200, Australia}

\email{Geoffrey.Campbell@anu.edu.au}

\thanks{Thanks to Professor Dr Henk Koppelaar, whose suggestions helped summarize, for this paper, parts of chapters 5, 12 and 21 of the author's draft book to appear in 2024.}

\keywords{Exact enumeration problems, generating functions. Partitions of integers. Elementary theory of partitions. Combinatorial identities, bijective combinatorics. Lattice points in specified regions.}
\subjclass[2010]{Primary: 05A15; Secondary: 05E40, 11Y11, 11P21}

\begin{abstract}
We set out an elementary approach to derive Visible Point Identities summed on lattice points of inverted triangle (2D), pyramid (3D), hyperpyramid (4D, 5D and so on) utilizing the greatest common divisor for the nD Visible Point Vectors. This enables study of partitions in nD space into vector parts distributed along straight lines radial from the origin in first hyperquadrant where coordinates of lattice points are all positive integers. We also give several new combinatorial identities for Visible Point Vector partitions.
\end{abstract}

\maketitle

\section{VPV identities in square hyperpyramid regions} \label{S:Intro VPV hyperpyramids}

In the decade up to 2000 the author published papers that culminated in the 2000 paper on hyperpyramid ``Visible Point Vector'' (VPV) identities. (See Campbell\index{Campbell, G.B.} \cite{gC1992,gC1993,gC1994a,gC1994b,gC1997,gC1998} then \cite{gC2000}.) These papers although breaking new ground, were not at the time taken any further than statement of a general theorem and a few prominent examples. Since then, these identities have not been developed further in the literature, despite there being evidently a large number of ways the Parametric Euler Sum Identities\index{Euler, L.} of the 21st century along with experimental computation results are definitely applicable. Partitions and lattice sums have been developed ``in their separate lanes'' in mathematical analysis, statistical mechanics physics, and in chemistry, without any great attempt to relate these areas of contemporary research. Add to this the possibility that light diffusion lattice models, random walk regimes, and stepping stone weighted partitions seem fundamentally applicable in contexts of VPV identities, and it becomes clear that the transition from integer partitions to vector partitions may be a path for future researches.

So, to be clear, we state the basic 2D, 3D and 4D identities as follows.

\begin{statement} For $|y|, |z|<1,$
  \begin{equation}\label{21.17i}
    \prod_{\substack{gcd(a,b)=1 \\ a<b \\ a \geq 0, \; b \geq 1}} \left( \frac{1}{1-y^a z^b} \right)^{\frac{1}{b}}
    = \left(\frac{1-yz}{1-z}\right)^{\frac{1}{1-y}}
  \end{equation}
        \begin{equation}  \nonumber
= 1 + \frac{z}{1!} + \begin{vmatrix}
    1 & -1 \\
    \frac{1-y^2}{1-y} & 1 \\
  \end{vmatrix} \frac{z^2}{2!}
  + \begin{vmatrix}
    1 & -1 & 0 \\
    \frac{1-y^2}{1-y} & 1 & -2 \\
    \frac{1-y^3}{1-y} & \frac{1-y^2}{1-y} & 1 \\
  \end{vmatrix} \frac{z^3}{3!}
+ \begin{vmatrix}
    1 & -1 & 0 & 0 \\
    \frac{1-y^2}{1-y} & 1 & -2 & 0 \\
    \frac{1-y^3}{1-y} & \frac{1-y^2}{1-y} & 1 & -3 \\
    \frac{1-y^4}{1-y} & \frac{1-y^3}{1-y} & \frac{1-y^2}{1-y} & 1 \\
  \end{vmatrix} \frac{z^4}{4!}
+ etc.
\end{equation}
   \end{statement}

\begin{statement} For each of $|x|, |y|, |z|<1,$
    \begin{equation}\label{21.18i}
    \prod_{\substack{gcd(a,b,c)=1 \\ a,b<c \\ a,b\geq0, \; c>0}} \left( \frac{1}{1-x^a y^b z^c} \right)^{\frac{1}{c}}
    = \left(\frac{(1-xz)(1-yz)}{(1-z)(1-xyz)}\right)^{\frac{1}{(1-x)(1-y)}}
  \end{equation}
      \begin{equation}  \nonumber
= 1 + \frac{z}{1!} + \begin{vmatrix}
    1 & -1 \\
    \frac{(1-x^2)(1-y^2)}{(1-x)(1-y)} & 1 \\
  \end{vmatrix} \frac{z^2}{2!}
  + \begin{vmatrix}
    1 & -1 & 0 \\
    \frac{(1-x^2)(1-y^2)}{(1-x)(1-y)} & 1 & -2 \\
    \frac{(1-x^3)(1-y^3)}{(1-x)(1-y)} & \frac{(1-x^2)(1-y^2)}{(1-x)(1-y)} & 1 \\
  \end{vmatrix} \frac{z^3}{3!}
      \end{equation}
  \begin{equation}  \nonumber
+ \begin{vmatrix}
    1 & -1 & 0 & 0 \\
    \frac{(1-x^2)(1-y^2)}{(1-x)(1-y)} & 1 & -2 & 0 \\
    \frac{(1-x^3)(1-y^3)}{(1-x)(1-y)} & \frac{(1-x^2)(1-y^2)}{(1-x)(1-y)} & 1 & -3 \\
    \frac{(1-x^4)(1-y^4)}{(1-x)(1-y)} & \frac{(1-x^3)(1-y^3)}{(1-x)(1-y)} & \frac{(1-x^2)(1-y^2)}{(1-x)(1-y)} & 1 \\
  \end{vmatrix} \frac{z^4}{4!}
+ etc.
\end{equation}
   \end{statement}

\begin{statement} For each of $|w|, |x|, |y|, |z|<1,$
    \begin{equation}\label{21.19i}
    \prod_{\substack{gcd(a,b,c,d)=1 \\ a,b,c<d \\ a,b,c\geq0, \; d>0}} \left( \frac{1}{1-w^a x^b y^c z^d} \right)^{\frac{1}{d}}
    = \left(\frac{(1-wz)(1-xz)(1-yz)(1-wxyz)}{(1-z)(1-wxz)(1-wyz)(1-xyz)}\right)^{\frac{1}{(1-w)(1-x)(1-y)}},
  \end{equation}
    \begin{equation}  \nonumber
= 1 + \frac{z}{1!} + \begin{vmatrix}
    1 & -1 \\
    \frac{(1-w^2)(1-x^2)(1-y^2)}{(1-w)(1-x)(1-y)} & 1 \\
  \end{vmatrix} \frac{z^2}{2!}
\end{equation}
  \begin{equation}  \nonumber
  + \begin{vmatrix}
    1 & -1 & 0 \\
    \frac{(1-w^2)(1-x^2)(1-y^2)}{(1-w)(1-x)(1-y)} & 1 & -2 \\
    \frac{(1-w^3)(1-x^3)(1-y^3)}{(1-w)(1-x)(1-y)} & \frac{(1-w^2)(1-x^2)(1-y^2)}{(1-w)(1-x)(1-y)} & 1 \\
  \end{vmatrix} \frac{z^3}{3!}
      \end{equation}
  \begin{equation}  \nonumber
+ \begin{vmatrix}
    1 & -1 & 0 & 0 \\
    \frac{(1-w^2)(1-x^2)(1-y^2)}{(1-w)(1-x)(1-y)} & 1 & -2 & 0 \\
    \frac{(1-w^3)(1-x^3)(1-y^3)}{(1-w)(1-x)(1-y)} & \frac{(1-w^2)(1-x^2)(1-y^2)}{(1-w)(1-x)(1-y)} & 1 & -3 \\
    \frac{(1-w^4)(1-x^4)(1-y^4)}{(1-w)(1-x)(1-y)} & \frac{(1-w^3)(1-x^3)(1-y^3)}{(1-w)(1-x)(1-y)} & \frac{(1-w^2)(1-x^2)(1-y^2)}{(1-w)(1-x)(1-y)} & 1 \\
  \end{vmatrix} \frac{z^4}{4!}
+ etc.
\end{equation}
   \end{statement}

Each of the above identities give us exact enumerations of certain functions of vector partitions. Suppose we say the $z$ variable is a ``vertical axis'' upon which to plot the 2D, 3D, 4D graph for the type of partitions and partition functions under consideration. Then the power series in $z$ give us exact determinant representations at each 1D, 2D or 3D layer corresponding to $z=1,2,3,...$ .

The proofs of the power series determinant coefficient functions in (\ref{21.17i}) to (\ref{21.19i}) rely only on applying Cramer's Rule to the coefficient recurrences, as well as differentiating the logarithms of both sides of the infinite products and their closed form evaluations.

So, in ensuing pages, we give the simplest $n$-space hyperpyramid VPV theorem due to the author in \cite{gC2000}. The so-called “Skewed Hyperpyramid $n$-space Identities" from \cite{gC2000} we shall cover in a later paper. The application of the determinant coefficient technique of our earlier work on hyperquadrant lattices is applicable and bearing some semblance to the $q$-binomial from earlier papers' variants. Note that for most identities in this paper, the left side products are taken over a set of integer lattice points inside an inverted 2D triangle lattice, or 3D pyramid, or 3+ dimension hyperpyramid on the Euclidean cartesian space.

In the first 15 years of the 21st century the summations found by the Borwein brothers Peter and Jonathan, their father David with their colleagues, see \cite{dB2006} to \cite{jB2013} have renewed interest in the old Euler Sums\index{Euler, L.}. Their results give us particular values of polylogarithms\index{Polylogarithm} and related functions popularized by Lewin \cite{lL1981,lL1982,lL1991}  involving the generalized Harmonic numbers\index{Harmonic number}. This work has been developed some way over nearly two decades so now we speak of the Mordell-Tornheim-Witten sums\index{Mordell-Tornheim-Witten (MTW) zeta function}, which are polylogarithm\index{Polylogarithm} generalizations all seen to be applicable to the VPV identities, but that connection is not yet fully worked through in the present literature. These newer researches including experimentally calculated results can, many of them, be substituted into VPV identities to give us exact results for weighted vector partitions. To make sense of these new results, we need to go back to fundamental definitions and ideas for partitions of vectors as distinct from those well considered already for integer partitions.

We examine the following correspondences. From the elementary generating function for unrestricted integer partitions we have,

\begin{equation} \label{21.20i}
  \prod_{n=1}^{\infty} \left(1+x^{1n}+x^{2n}+x^{3n}+ \ldots \right)=1+p(1)x+p(2)x^2+p(3)x^3+\ldots .
\end{equation}

So, equations (\ref{21.17i}) to (\ref{21.19i}) similarly imply

\begin{equation} \label{21.21i}
  \prod_{\substack{gcd(a,b)=1; \; a<b \\ a \geq 0, \; b \geq 1 }}
  \left( 1+\binom{1/b}{1}(y^a z^b)+\binom{1/b}{2}(y^a z^b)^2+\binom{1/b}{3}(y^a z^b)^3+\ldots \right)
\end{equation}
\begin{equation*}
  = \sum_{\substack{a<b \\ a \geq 0, \; b \geq 1 }} V_2(a,b) y^a z^b.
\end{equation*}

\begin{equation} \label{21.22i}
  \prod_{\substack{gcd(a,b,c)=1; \\ a,b<c \\ a,b \geq 0, \; c \geq 1 }}
  \left( 1+\binom{1/c}{1}(x^a y^b z^c)+\binom{1/c}{2}(x^a y^b z^c)^2+\binom{1/c}{3}(x^a y^b z^c)^3+\ldots \right)
\end{equation}
\begin{equation*}
  = \sum_{\substack{a,b<c \\ a,b \geq 0, \; c \geq 1 }} V_3(a,b,c) x^a y^b z^c.
\end{equation*}

\begin{equation} \label{21.23i}
  \prod_{\substack{gcd(a,b,c,d)=1; \\ a,b,c<d \\ a,b,c \geq 0, \; d \geq 1 }}
  \left( 1+\binom{1/d}{1}(w^a x^b y^c z^d)+\binom{1/d}{2}(w^a x^b y^c z^d)^2+\binom{1/d}{3}(w^a x^b y^c z^d)^3+\ldots \right)
\end{equation}
\begin{equation*}
  = \sum_{\substack{a,b,c<d \\ a,b,c \geq 0, \; d \geq 1 }} V_4(a,b,c,d) w^a x^b y^c z^d.
\end{equation*}

Equations (\ref{21.21i}) to (\ref{21.23i}) ought to give us exact closed form generating functions $V_2$, $V_3$, $V_4$, in 2D, 3D and 4D space respectively. We put this aside for now.

\newpage

\section{Vector Partitions whose parts are on \MakeLowercase{n}D straight lines.}

From equation (\ref{21.20i}) we see that replacing $x$ by $y^a z^b$ for $|y^a z^b|<1$ where $a$ and $b$ are coprime positive integers, gives the equation
\begin{equation} \label{21.24i}
  \prod_{n=1}^{\infty} \frac{1}{1 - (y^a z^b)^n}
  =\prod_{n=1}^{\infty} \left(1+(y^a z^b)^{1n}+(y^a z^b)^{2n}+(y^a z^b)^{3n}+ \ldots \right)
 \end{equation}
\begin{equation*}
  =1+p(1)(y^a z^b)^1+p(2)(y^a z^b)^2+p(3)(y^a z^b)^3+\ldots .
\end{equation*}

Of course this is just a thinly disguised version of the generating function for unrestricted integer partitions. ie. one-dimensional partitions. In the 2D case we are saying that the number of partitions of an integer lattice point vector $\langle A,B \rangle$ on the line $z=ay/b$ for $gcd(a,b)=1$ into vector parts also on this line, is equal to $p(gcd(A,B))$.

So (\ref{21.24i}), enumerating the number of vector partitions for the lattice points along a 1D line in 2D space, applies equally to a 1D line in 3D space as follows, where $a$, $b$ and $c$ are positive integers with $gcd(a,b,c)=1$.
\begin{equation} \label{21.25i}
  \prod_{n=1}^{\infty} \frac{1}{1 - (x^a y^b z^c)^n}
  =\prod_{n=1}^{\infty} \left(1+(x^a y^b z^c)^{1n}+(x^a y^b z^c)^{2n}+(x^a y^b z^c)^{3n}+ \ldots \right)
 \end{equation}
\begin{equation*}
  =1+p(1)(x^a y^b z^c)^1+p(2)(x^a y^b z^c)^2+p(3)(x^a y^b z^c)^3+\ldots .
\end{equation*}

Similarly, in the 4D case extension corresponding to the 3D equation (\ref{21.25i}) we have, where $a$, $b$, $c$ and $d$ are positive integers with gcd(a,b,c,d)=1,
\begin{equation} \label{21.26i}
  \prod_{n=1}^{\infty} \frac{1}{1 - (w^a x^b y^c z^d)^n}
   \end{equation}
\begin{equation*}
  =\prod_{n=1}^{\infty} \left(1+(w^a x^b y^c z^d)^{1n}+(w^a x^b y^c z^d)^{2n}+(w^a x^b y^c z^d)^{3n}+ \ldots \right)
 \end{equation*}
\begin{equation*}
  =1+p(1)(w^a x^b y^c z^d)^1+p(2)(w^a x^b y^c z^d)^2+p(3)(w^a x^b y^c z^d)^3+\ldots .
\end{equation*}

We shall now do something interesting with equation (\ref{21.24i}). We write an ``Upper VPV'' identity derived from it as follows.

\begin{equation} \label{21.27i}
  \prod_{\substack{a,b,n \geq 1; \; a \leq b \\ gcd(a,b)=1}} \frac{1}{1 - (y^a z^b)^n}
 \end{equation}
 \begin{eqnarray*}
     &=& \frac{1}{1 - y^1 z^1} \\
     &\times& \frac{1}{1 - y^1 z^2} \frac{1}{1 - y^2 z^2} \\
     &\times& \frac{1}{1 - y^1 z^3} \frac{1}{1 - y^2 z^3} \frac{1}{1 - y^3 z^3} \\
     &\times& \frac{1}{1 - y^1 z^4} \frac{1}{1 - y^2 z^4} \frac{1}{1 - y^3 z^4} \frac{1}{1 - y^4 z^4} \\
     &\times& \frac{1}{1 - y^1 z^5} \frac{1}{1 - y^2 z^5} \frac{1}{1 - y^3 z^5} \frac{1}{1 - y^4 z^5} \frac{1}{1 - y^5 z^5} \\
     &\times& \frac{1}{1 - y^1 z^6} \frac{1}{1 - y^2 z^6} \frac{1}{1 - y^3 z^6} \frac{1}{1 - y^4 z^6} \frac{1}{1 - y^5 z^6} \frac{1}{1 - y^6 z^6} \\
     &\times& etc.
 \end{eqnarray*}
\begin{equation*}
  =\prod_{\substack{a,b,n \geq 1; \; a \leq b \\ gcd(a,b)=1}} \left(1+(y^a z^b)^{1n}+(y^a z^b)^{2n}+(y^a z^b)^{3n}+ \ldots \right)
 \end{equation*}
\begin{equation*}
  =\prod_{\substack{a,b \geq 1; \; a \leq b \\ gcd(a,b)=1}} \left( 1+p(1)(y^a z^b)^1+p(2)(y^a z^b)^2+p(3)(y^a z^b)^3+\ldots \right)
\end{equation*}
\begin{equation*}
  = \sum_{n_1,n_2,n_3... \geq 0} p(n_1)(y^1 z^1)^{n_1} p(n_2)(y^1 z^2)^{n_2}p(n_3)(y^1 z^3)^{n_3}p(n_4)(y^2 z^3)^{n_3}\ldots
\end{equation*}
\begin{equation*}
  = \sum_{n_1,n_2,n_3... \geq 0} p(n_1)p(n_2)p(n_3)p(n_4)\cdots(y^1 z^1)^{n_1}(y^1 z^2)^{n_2}(y^1 z^3)^{n_3}(y^2 z^3)^{n_3}\ldots
\end{equation*}
\begin{equation*}
  = \sum_{n_1,n_2,n_3... \geq 0} p(n_1)p(n_2)p(n_3)p(n_4)\cdots(y^{(1n_1+1n_2+1n_3+2n_4+\ldots)}z^{(1n_1+2n_2+3n_3+3n_4+\ldots)})
\end{equation*}
where each coefficient of $n_i$ in the index sum of $y$ is coprime to the coefficient of $n_i$ in the index sum of $z$;
 \begin{eqnarray*}
     &=& 1 + p_{(1|1)}y^1 z^1 \\
     & & + \; p_{(1|2)} y^1 z^2 + p_{(2|2)} y^2 z^2 \\
     & & + \; p_{(1|3)} y^1 z^3 + p_{(2|3)} y^2 z^3 + p_{(3|3)} y^3 z^3 \\
     & & + \; p_{(1|4)} y^1 z^4 + p_{(2|4)} y^2 z^4 + p_{(3|4)} y^3 z^4 + p_{(4|4)} y^4 z^4 \\
     & & + \; p_{(1|5)} y^1 z^5 + p_{(2|5)} y^2 z^5 + p_{(3|5)} y^3 z^5 + p_{(4|5)} y^4 z^5 + p_{(5|5)} y^5 z^5 \\
     & & + \; etc.
 \end{eqnarray*}
where $p_{(a|b)} := p_{(a|b)}(y,z)$ is the coefficient function of $y^a z^b$.

This enables study of partitions in $n$D space into vector parts distributed along straight lines radial from the origin in first hyperquadrant where coordinates of lattice points are all positive integers.

\bigskip

\textbf{2D Vector Partitions whose parts are on two straight lines.}

The above analysis is with respect to vector partitions in the upper half first quadrant above and including the line $y=z$. The simplest version of this theory departing from the single radial from the origin line of lattice points, would be to state the result if dealing with 2D partitions from two such radial lines of lattice points. For example, consider the two lines with equations $y=z/2$ and $y=z/3$. The lattice point vectors along these lines in the first quadrant may be listed as
\begin{eqnarray*}
  S_1 &=& \{\langle 1,2\rangle ,\langle 2,4\rangle ,\langle 3,6\rangle ,\langle 4,8\rangle ,\langle 5,10\rangle ,\ldots\}; \\
  S_2 &=& \{\langle 1,3\rangle ,\langle 2,6\rangle ,\langle 3,9\rangle ,\langle 4,12\rangle ,\langle 5,15\rangle ,\ldots\}.
\end{eqnarray*}
Following the above rationale we see that the generating function for 2D vector partitions into parts contained in $S_1$ and $S_2$ is
\begin{equation*}
  \frac{1}{((1-yz^2)(1-y^2z^4)(1-y^3z^6)\cdots)(1-yz^3)(1-y^2z^6)(1-y^3z^9)\cdots)}
\end{equation*}
\begin{equation*}
 = \left(1+p(1)y^1 z^2+p(2)y^2 z^4+p(3)y^3 z^6+\ldots \right)\left(1+p(1)y^1 z^3+p(2)y^2 z^6+p(3)y^3 z^9+\ldots \right)
\end{equation*}
\begin{eqnarray*}
    &=&  1 + p(1) y z^2 + p(1) y z^3 \\
    & &  + \; p(2) y^2 z^4 + p(1)^2 y^2 z^5 + p(2) y^2 z^6 \\
    & &  + \; p(3) y^3 z^6 + p(1) p(2) y^3 z^7 + p(1) p(2) y^3 z^8 + p(3) y^3 z^9 \\
    & &  + \; etc.
\end{eqnarray*}
\begin{eqnarray*}
    &=&  1 + y z^2 + y z^3 + 2 y^2 z^4 + y^2 z^5 + y^2 (3 y + 2) z^6 \\
    & &  + \; 2 y^3 z^7 + y^3 (5 y + 2) z^8 + 3 y^3 (y + 1) z^9 + y^4 (7 y + 4) z^{10} \\
    & &  + \; y^4 (5 y + 3) z^{11} + y^4 (11 y^2 + 6 y + 5) z^{12} + y^5 (7 y + 6) z^{13} \\
    & &  + \; y^5 (14 y^2 + 10 y + 5) z^{14} + y^5 (11 y^2 + 9 y + 7) z^{15} \\
    & &  + \; etc.
\end{eqnarray*}
The coefficients here plot onto the grid
\begin{equation*}
\begin{array}{c|ccccccccc}
\vdots & \vdots & \vdots & \vdots & \vdots & \vdots & \vdots & \vdots & \vdots & \vdots  \\
15 &   &   &   &   &   & 7 & 9 & 11 &   \cdots \\
14 &   &   &   &   &   & 5 & 10 & 14 &   \cdots \\
13 &   &   &   &   &   & 6 & 7 &   &   \cdots \\
12 &   &   &   &   & 5 & 6 & 11 &   &   \cdots \\
11 &   &   &   &   & 3 & 5 &   &   &   \cdots \\
10 &   &   &   &   & 4 & 7 &   &   &   \cdots \\
9  &   &   &   & 3 & 3 &   &   &   &   \cdots \\
8  &   &   &   & 2 & 5 &   &   &   &   \cdots \\
7  &   &   &   & 2 &   &   &   &   &   \cdots \\
6  &   &   & 2 & 3 &   &   &   &   &   \cdots \\
5  &   &   & 1 &   &   &   &   &   &   \cdots \\
4  &   &   & 2 &   &   &   &   &   &   \cdots \\
3  &   & 1 &   &   &   &   &   &   &   \cdots \\
2  &   & 1 &   &   &   &   &   &   &   \cdots \\
1  &   &   &   &   &   &   &   &   &   \cdots \\
0  & 1 &   &   &   &   &   &   &   &   \cdots \\ \hline
z/y& 0 & 1 & 2 & 3 & 4 & 5 & 6 & 7 &  \cdots
  \end{array}
\end{equation*}

\textbf{Example interpretations reading from this graph.} 
\begin{enumerate}
  \item [1)] The number of partitions of vector $\langle 7,15 \rangle$ into unrestricted number of parts from $S_1$ and $S_2$ is 11.
  \item [2)] The number of partitions of vector $\langle 5,10 \rangle$ into unrestricted number of parts from $S_1$ and $S_2$ is 7.
  \item [3)] The number of partitions of vector $\langle 4,9 \rangle$ into unrestricted number of parts from $S_1$ and $S_2$ is 3.
  \item [4)] The number of partitions of vector $\langle 4,7 \rangle$ into unrestricted number of parts from $S_1$ and $S_2$ is 0.
\end{enumerate}

\bigskip

\textbf{3D Vector Partitions whose parts are on two straight lines.}

A further simple version of the above approach allows us to work out exactly the number of partitions into vector parts that lie upon ``radial from the origin'' lines of lattice points. We state the following example in 3D partitions in the first 3D hyperquadrant (ie. with lattice points whose co-ordinates are triples of positive integers.) For example, consider the first line through 3D space with equation $x=y/2=z/3$; then the second line with equation $x=y/3=z/4$. The lattice point vectors along these lines in the first quadrant may be listed as
\begin{eqnarray*}
  S_3 &=& \{\langle 1,2,3\rangle ,\langle 2,4,6\rangle ,\langle 3,6,9\rangle ,\langle 4,8,12\rangle ,\langle 5,10,15\rangle ,\ldots\}; \\
  S_4 &=& \{\langle 1,3,4\rangle ,\langle 2,6,8\rangle ,\langle 3,9,12\rangle ,\langle 4,12,16\rangle ,\langle 5,15,20\rangle ,\ldots\}.
\end{eqnarray*}
Applying our rationale we see that the generating function for 3D vector partitions into parts contained in $S_3$ and $S_4$ is
\begin{equation*}
  \frac{1}{((1-xy^2z^3)(1-x^2y^4z^6)(1-x^3y^6z^9)\cdots)(1-xy^3z^4)(1-x^2y^6z^8)(1-x^3y^9z^{12})\cdots)}
\end{equation*}
\begin{equation*}
 = \left(1+p(1)xy^2z^3+p(2)x^2y^4z^6+p(3)x^3y^6z^9+\ldots \right) \quad \quad \quad
 \end{equation*}
\begin{equation*}
  \quad \quad \times \left(1+p(1)xy^3z^4+p(2)x^2y^6z^8+p(3)x^3y^9z^{12}+\ldots \right)
\end{equation*}
\begin{eqnarray*}
    &=&  1 + p(1) x y^3 z^4+ p(1) x y^2 z^3 \\
    & &  + \; p(2) x^2 y^6 z^8 + p(1)^2 x^2 y^5 z^7 + p(2) x^2 y^4 z^6  \\
    & &  + \; p(3) x^3 y^6 z^9 + p(1) p(2) x^3 y^7 z^{10} + p(1) p(2) x^3 y^8 z^11 + p(3) x^3 y^9 z^{12} \\
    & &  + \; etc.
\end{eqnarray*}

\bigskip

There are many extended possibilities for the above enumerations of partitions in higher cartesian space sets of lattice point vectors.  There is no reason other than a formal complexity why we can't create exact enumerative formulas for ``$n$D Vector Partitions whose parts are on $m$ straight lines'' for $m$ and $n$ both arbitrary fixed positive integers.

\bigskip

\textbf{Distinct Vector Partitions along an \textit{n}-space line.}

Recall that Euler, in addition to giving us the ``unrestricted'' integer partitions generating function, also noted that for $|x|<1$,

\begin{equation}\label{21.28i}
  \prod_{k=1}^{\infty} (1+x^k) = \sum_{n=1}^{\infty} \mathcal{D}(n)x^n,
\end{equation}

where $\mathcal{D}(n)$ is the number of partitions of positive integer $n$ into distinct integer parts.

From equation (\ref{21.28i}) we see that replacing $x$ by $y^a z^b$ for $|y^a z^b|<1$ where $a$ and $b$ are coprime positive integers, gives the equation
\begin{equation} \label{21.29i}
  \prod_{n=1}^{\infty} (1 + (y^a z^b)^n)
  =\prod_{n=1}^{\infty} \left(1+(y^a z^b)^{1n}+(y^a z^b)^{2n}+(y^a z^b)^{3n}+ \ldots \right)
 \end{equation}
\begin{equation*}
  =1+\mathcal{D}(1)(y^a z^b)^1+\mathcal{D}(2)(y^a z^b)^2+\mathcal{D}(3)(y^a z^b)^3+\ldots .
\end{equation*}

Similarly to the previous rationale using ``unrestricted'' partitions, this is a version of the generating function for integer partitions into distinct vector parts. ie. one-dimensional partitions into distinct parts. In the 2D case we are saying that the number of partitions of an integer lattice point vector $\langle A,B \rangle$ on the line $z=ay/b$ for $gcd(a,b)=1$ into ``distinct'' vector parts also on this line, is equal to $\mathcal{D}(gcd(A,B))$.

In a further example, consider again the two lines with equations $y=z/2$ and $y=z/3$. The lattice point vectors along these lines in the first quadrant we list again as
\begin{eqnarray*}
  S_1 &=& \{\langle 1,2\rangle ,\langle 2,4\rangle ,\langle 3,6\rangle ,\langle 4,8\rangle ,\langle 5,10\rangle ,\ldots\}; \\
  S_2 &=& \{\langle 1,3\rangle ,\langle 2,6\rangle ,\langle 3,9\rangle ,\langle 4,12\rangle ,\langle 5,15\rangle ,\ldots\}.
\end{eqnarray*}
Following the above rationale we see that the generating function for 2D vector partitions into distinct parts contained in $S_1$ and $S_2$ is
\begin{equation*}
  ((1+yz^2)(1+y^2z^4)(1+y^3z^6)\cdots)((1+yz^3)(1+y^2z^6)(1+y^3z^9)\cdots)
\end{equation*}
\begin{equation*}
 = \left(1+\mathcal{D}(1)y^1 z^2+\mathcal{D}(2)y^2 z^4+\mathcal{D}(3)y^3 z^6+\ldots \right) \quad \quad
 \end{equation*}
\begin{equation*}
 \quad \times \left(1+\mathcal{D}(1)y^1 z^3+\mathcal{D}(2)y^2 z^6+\mathcal{D}(3)y^3 z^9+\ldots \right)
\end{equation*}
\begin{eqnarray*}
    &=&  1 + \mathcal{D}(1) y z^2 + \mathcal{D}(1) y z^3 \\
    & &  + \; \mathcal{D}(2) y^2 z^4 + \mathcal{D}(1)^2 y^2 z^5 + \mathcal{D}(2) y^2 z^6 \\
    & &  + \; \mathcal{D}(3) y^3 z^6 + \mathcal{D}(1) \mathcal{D}(2) y^3 z^7 + \mathcal{D}(1) \mathcal{D}(2) y^3 z^8 + \mathcal{D}(3) y^3 z^9 \\
    & &  + \; etc.
\end{eqnarray*}
\begin{eqnarray*}
    &=&  1 + y z^2 + y z^3 + y^2 z^4 + y^2 z^5 + y^2 (2 y + 1) z^6 \\
    & &  + \; y^3 z^7 + y^3 (2 y + 1) z^8 + 2 y^3 (y + 1) z^9 + y^4 (3 y + 1) z^{10} \\
    & &  + \; 2 y^4 (y + 1) z^{11} + 2 y^4 (2 y^2 + y + 1) z^{12} + y^5 (3 y + 2) z^13 \\
    & &  + \; 2 y^5 (2 y^2 + y + 1) z^{14} + y^5 (4 y^2 + 4 y + 3) z^{15} \\
    & &  + \; etc.
\end{eqnarray*}
We can easily plot the coefficients here onto the grid
\begin{equation*}
\begin{array}{c|ccccccccc}
\vdots & \vdots & \vdots & \vdots & \vdots & \vdots & \vdots & \vdots & \vdots & \vdots  \\
15 &   &   &   &   &   & 3 & 4 & 4 &   \cdots \\
14 &   &   &   &   &   & 2 & 2 & 4 &   \cdots \\
13 &   &   &   &   &   & 2 & 3 &   &   \cdots \\
12 &   &   &   &   & 2 & 2 & 4 &   &   \cdots \\
11 &   &   &   &   & 2 & 2 &   &   &   \cdots \\
10 &   &   &   &   & 1 & 3 &   &   &   \cdots \\
9  &   &   &   & 2 & 2 &   &   &   &   \cdots \\
8  &   &   &   & 1 & 2 &   &   &   &   \cdots \\
7  &   &   &   & 1 &   &   &   &   &   \cdots \\
6  &   &   & 1 & 2 &   &   &   &   &   \cdots \\
5  &   &   & 1 &   &   &   &   &   &   \cdots \\
4  &   &   & 1 &   &   &   &   &   &   \cdots \\
3  &   & 1 &   &   &   &   &   &   &   \cdots \\
2  &   & 1 &   &   &   &   &   &   &   \cdots \\
1  &   &   &   &   &   &   &   &   &   \cdots \\
0  & 1 &   &   &   &   &   &   &   &   \cdots \\ \hline
z/y& 0 & 1 & 2 & 3 & 4 & 5 & 6 & 7 &  \cdots
  \end{array}
\end{equation*}
\textbf{Example interpretations reading from this graph.}
\begin{enumerate}
  \item [1)] The number of partitions of vector $\langle 7,15 \rangle$ using distinct parts from $S_1$ and $S_2$ is 11.
  \item [2)] The number of partitions of vector $\langle 5,10 \rangle$ using distinct parts from $S_1$ and $S_2$ is 7.
  \item [3)] The number of partitions of vector $\langle 4,9 \rangle$ using distinct parts from $S_1$ and $S_2$ is 3.
  \item [4)] The number of partitions of vector $\langle 4,7 \rangle$ using distinct parts from $S_1$ and $S_2$ is 0.
\end{enumerate}

  \bigskip

\bigskip

\section{Deriving 2D VPV identities in extended triangle regions} \label{S:2D VPV hyperpyramids}

In this section we derive the $2D$ Visible Point Vector identities from essentially creating them out of a simple summation transformation based on the simple idea that the integer lattice points in the first quadrant have co-ordinates that are either coprime integer pairs namely, ``lattice points visible from the origin'', or co-ordinates that are the integer multiples of the coprime pairs. As we did in the hyperquadrant paper, we again start with a simple $2D$ summation. Consider

\begin{equation}  \nonumber
   \sum_{n=1}^{\infty}  \left( \sum_{m=1}^{n} \frac{y^m}{m^a} \right) \frac{z^n}{n^b}
\end{equation}

\begin{equation}  \nonumber
=\left(\frac{y^1}{1^a}\right)\frac{z^1}{1^b}
+\left(\frac{y^1}{1^a}+\frac{y^2}{2^a}\right)\frac{z^2}{2^b}
+\left(\frac{y^1}{1^a}+\frac{y^2}{2^a}+\frac{y^3}{3^a}\right)\frac{z^3}{3^b}
+\left(\frac{y^1}{1^a}+\frac{y^2}{2^a}+\frac{y^3}{3^a}+\frac{y^4}{4^a}\right)\frac{z^4}{4^b}
\end{equation}
\begin{equation}  \nonumber
+\left(\frac{y^1}{1^a}+\frac{y^2}{2^a}+\frac{y^3}{3^a}+\frac{y^4}{4^a}+\frac{y^5}{5^a}\right)\frac{z^5}{5^b}
+\left(\frac{y^1}{1^a}+\frac{y^2}{2^a}+\frac{y^3}{3^a}+\frac{y^4}{4^a}+\frac{y^5}{5^a}+\frac{y^6}{6^a}\right)\frac{z^6}{6^b}+\cdots
\end{equation}

\bigskip

\begin{equation}  \nonumber
 =\frac{y^1 z^1}{1^a 1^b}
\end{equation}
\begin{equation}  \nonumber
 +\frac{y^1 z^2}{1^a 2^b}+\frac{y^2 z^2}{2^a 2^b}
\end{equation}
\begin{equation}  \nonumber
 +\frac{y^1 z^3}{1^a 3^b}+\frac{y^2 z^3}{2^a 3^b}+\frac{y^3 z^3}{3^a 3^b}
\end{equation}
\begin{equation}  \nonumber
 +\frac{y^1 z^4}{1^a 4^b}+\frac{y^2 z^4}{2^a 4^b}+\frac{y^3 z^4}{3^a 4^b}+\frac{y^4 z^4}{4^a 4^b}
\end{equation}
\begin{equation}  \nonumber
 +\frac{y^1 z^5}{1^a 5^b}+\frac{y^2 z^5}{2^a 5^b}+\frac{y^3 z^5}{3^a 5^b}+\frac{y^4 z^5}{4^a 5^b}+\frac{y^5 z^5}{5^a 5^b}
\end{equation}
\begin{equation}  \nonumber
 +\frac{y^1 z^6}{1^a 6^b}+\frac{y^2 z^6}{2^a 6^b}+\frac{y^3 z^6}{3^a 6^b}+\frac{y^4 z^6}{4^a 6^b}+\frac{y^5 z^6}{5^a 6^b}+\frac{y^6 z^6}{6^a 6^b}
\end{equation}
\begin{equation}  \nonumber
 +\frac{y^1 z^7}{1^a 7^b}+\frac{y^2 z^7}{2^a 7^b}+\frac{y^3 z^7}{3^a 7^b}+\frac{y^4 z^7}{4^a 7^b}+\frac{y^5 z^7}{5^a 7^b}+\frac{y^6 z^7}{6^a 7^b}+\frac{y^7 z^7}{7^a 7^b}
\end{equation}
\begin{equation}  \nonumber
 + \quad \vdots \quad + \quad \vdots \quad + \quad \vdots \quad + \quad \vdots \quad + \quad \vdots \quad + \quad \vdots \quad + \quad \vdots \quad   \ddots
\end{equation}
\begin{equation}  \nonumber
 = \sum_{m,n \geq 1; m \leq n}^{\infty}  \frac{y^m z^n}{m^a n^b}
\end{equation}
\begin{equation}  \nonumber
 = \sum_{\substack{ h,j,k \geq 1 \\ j \leq k ; \, (j,k)=1}}  \frac{(y^j z^k)^h}{h^{a+b} (j^a k^b)}
\end{equation}
\begin{equation}  \nonumber
 = \sum_{\substack{ j,k \geq 1 \\ j \leq k ; \, (j,k)=1}}  \frac{1}{(j^a k^b)}   \sum_{h=1}^{\infty} \frac{(y^j z^k)^h}{h^{a+b}}
\end{equation}
\begin{equation}  \nonumber
 = \sum_{\substack{ j,k \geq 1 \\ j \leq k ; \, (j,k)=1}} \frac{1}{(j^a k^b)}   \log \left( \frac{1}{1 - y^j z^k} \right) \quad if \quad a+b=1.
\end{equation}

Therefore, we have shown that
\begin{equation}  \nonumber
 \sum_{n=1}^{\infty}  \left( \sum_{m=1}^{n} \frac{y^m}{m^a} \right) \frac{z^n}{n^b} = \sum_{\substack{ j,k \geq 1 \\ j \leq k ; \, (j,k)=1}} \frac{1}{(j^a k^b)}   \log \left( \frac{1}{1 - y^j z^k} \right) \quad if \quad a+b=1.
\end{equation}

 Exponentiating both sides gives us the $2D$ first extended triangle VPV identity, where in this $2D$ case the $nD$ pyramid reduces to the form of a triangle shaped array of lattice point vectors, and so we can state the

 \begin{theorem}     \label{vpv-pyramid2D-thm}
 \textbf{The $2D$ first quadrant triangle VPV identity.} For $|y|<1, |z|<1,$
 \begin{equation}   \label{21.01}
    \prod_{\substack{ j,k \geq 1 \\ j \leq k ; \, (j,k)=1}} \left( \frac{1}{1-y^j z^k} \right)^{\frac{1}{j^a k^b}}
    = \exp\left\{ \sum_{n=1}^{\infty}  \left( \sum_{m=1}^{n} \frac{y^m}{m^a} \right) \frac{z^n}{n^b} \right\} \quad if \quad a+b=1.
  \end{equation}
 \end{theorem}

As with our earlier exploits into the $2D$ first quadrant case, for the present result we take some simple example cases where new and interesting results arise.

So, let us take the case where $a=0, b=1$, giving us

 \begin{equation}   \nonumber
    \prod_{\substack{ j,k \geq 1 \\ j \leq k ; \, (j,k)=1}} \left( \frac{1}{1-y^j z^k} \right)^{\frac{1}{k}}
    = \exp\left\{ \sum_{n=1}^{\infty}  \left( \sum_{m=1}^{n} y^m \right) \frac{z^n}{n} \right\}
  \end{equation}
 \begin{equation}   \nonumber
    = \exp\left\{ \sum_{n=1}^{\infty}  \left( y \frac{1-y^n}{1-y} \right) \frac{z^n}{n} \right\}
    = \exp\left\{ \frac{y}{1-y} \log \left( \frac{1-yz}{1-z} \right)   \right\}.
  \end{equation}

So, we arrive then at the following pair of equivalent results,

 \begin{equation}   \label{21.02}
    \prod_{\substack{ j,k \geq 1 \\ j \leq k ; \, (j,k)=1}} \left( \frac{1}{1-y^j z^k} \right)^{\frac{1}{k}}
        =  \left( \frac{1-yz}{1-z} \right)^{\frac{y}{1-y}} ,
  \end{equation}

 and

 \begin{equation}   \label{21.03}
    \prod_{\substack{ j,k \geq 1 \\ j \leq k ; \, (j,k)=1}} \left( 1-y^j z^k \right)^{\frac{1}{k}}
        =  \left( \frac{1-z}{1-yz} \right)^{\frac{y}{1-y}} .
  \end{equation}

From here, multiply both sides of (\ref{21.02}) and the case of (\ref{21.03}) with $y \mapsto y^2$ and $z \mapsto z^2$ to get,

 \begin{equation}   \label{21.04}
    \prod_{\substack{ j,k \geq 1; \, j \leq k \\ gcd(j,k)=1}} \left( 1+y^j z^k \right)^{\frac{1}{k}}
        =  \left( \frac{1-yz}{1-z} \right)^{\frac{y}{1-y}} \left( \frac{1-z^2}{1-y^2z^2} \right)^{\frac{y^2}{1-y^2}} .
  \end{equation}

Obviously multiplying both sides of (\ref{21.03}) and (\ref{21.04}) give a restated (\ref{21.04}).

Particular cases: $y = \frac{1}{2}$ gives us from (\ref{21.03}) and (\ref{21.04}) the remarkable two results,

 \begin{equation}   \nonumber
    \prod_{\substack{ j,k \geq 1; \, j \leq k \\ gcd(j,k)=1}} \left( 1- \frac{z^k}{2^j} \right)^{\frac{1}{k}}
        =   \frac{2-2z}{2-z}  = 1 - \frac{z}{2} - \frac{z^2}{4} - \frac{z^3}{8} - \frac{z^4}{16} - \frac{z^5}{32} - \ldots
  \end{equation}
 \begin{equation}   \nonumber
 = \left( 1- \frac{z}{2} \right)
  \end{equation}
 \begin{equation}   \nonumber
 \sqrt{\left( 1- \frac{z^2}{2^1} \right)}
  \end{equation}
 \begin{equation}   \nonumber
 \sqrt[3]{\left( 1- \frac{z^3}{2^1} \right)\left( 1- \frac{z^3}{2^2} \right)}
  \end{equation}
 \begin{equation}   \nonumber
 \sqrt[4]{\left( 1- \frac{z^4}{2^1} \right)\left( 1- \frac{z^4}{2^3} \right)}
  \end{equation}
 \begin{equation}   \nonumber
 \sqrt[5]{\left( 1- \frac{z^5}{2^1} \right)\left( 1- \frac{z^5}{2^2} \right)\left( 1- \frac{z^5}{2^3} \right)\left( 1- \frac{z^5}{2^4} \right)}
  \end{equation}
 \begin{equation}   \nonumber
 \sqrt[6]{\left( 1- \frac{z^6}{2^1} \right)\left( 1- \frac{z^6}{2^5} \right)}
 \end{equation}
 \begin{equation}   \nonumber
 \vdots \, ,
 \end{equation}

 \begin{equation}   \nonumber
    \prod_{\substack{ j,k \geq 1; \, j \leq k \\ gcd(j,k)=1}} \left( 1+ \frac{z^k}{2^j} \right)^{\frac{1}{k}}
    = \frac{1-\frac{z}{2}}{1-z}\sqrt[3]{\frac{1-z^2}{1-\frac{z^2}{4}}} = 1 +\frac{z}{2} +\frac{z^2}{4} +\frac{3z^3}{8} +\frac{z^4}{4} +\frac{5z^5}{16} + \ldots
  \end{equation}
 \begin{equation}   \nonumber
 = \left( 1+ \frac{z}{2} \right)
  \end{equation}
 \begin{equation}   \nonumber
 \sqrt{\left( 1+ \frac{z^2}{2^1} \right)}
  \end{equation}
 \begin{equation}   \nonumber
 \sqrt[3]{\left( 1+ \frac{z^3}{2^1} \right)\left( 1+ \frac{z^3}{2^2} \right)}
  \end{equation}
 \begin{equation}   \nonumber
 \sqrt[4]{\left( 1+ \frac{z^4}{2^1} \right)\left( 1+ \frac{z^4}{2^3} \right)}
  \end{equation}
 \begin{equation}   \nonumber
 \sqrt[5]{\left( 1+ \frac{z^5}{2^1} \right)\left( 1+ \frac{z^5}{2^2} \right)\left( 1+ \frac{z^5}{2^3} \right)\left( 1+ \frac{z^5}{2^4} \right)}
  \end{equation}
 \begin{equation}   \nonumber
 \sqrt[6]{\left( 1+ \frac{z^6}{2^1} \right)\left( 1+ \frac{z^6}{2^5} \right)}
 \end{equation}
 \begin{equation}   \nonumber
 \vdots .
 \end{equation}

These two equations can be easily verified on a calculating engine like Mathematica or WolframAlpha by expanding each side into it's Taylor series around $z=0$ and comparing coefficients of like powers of $z$. Next, take the cases of (\ref{21.03}) and (\ref{21.04}) with $y=2$, both of which converge if $|z|<\frac{1}{2}$, so then, after a slight adjustment to both sides we have,

 \begin{equation}   \nonumber
    \prod_{\substack{ j,k \geq 1; \, j < k \\ gcd(j,k)=1}} \left( 1- 2^j z^k  \right)^{\frac{1}{k}}
        =  1- \frac{z^2}{(1-z)^2} = 1 - z^2 - 2z^3 - 3z^4 - 4z^5 - \ldots - n z^{n+1} - \ldots
  \end{equation}
 \begin{equation}   \nonumber
 = \sqrt{\left( 1- 2^1 z^2 \right)}
  \end{equation}
 \begin{equation}   \nonumber
 \sqrt[3]{\left( 1- 2^1 z^3 \right)\left( 1- 2^2 z^3 \right)}
  \end{equation}
 \begin{equation}   \nonumber
 \sqrt[4]{\left( 1- 2^1 z^4 \right)\left( 1- 2^3 z^4 \right)}
  \end{equation}
 \begin{equation}   \nonumber
 \sqrt[5]{\left( 1- 2^1 z^5 \right)\left( 1- 2^2 z^5 \right)\left( 1- 2^3 z^5 \right)\left( 1- 2^4 z^5 \right)}
  \end{equation}
 \begin{equation}   \nonumber
 \sqrt[6]{\left( 1- 2^1 z^6 \right)\left( 1- 2^5 z^6 \right)}
 \end{equation}
 \begin{equation}   \nonumber
 \sqrt[7]{\left( 1- 2^1 z^7 \right)\left( 1- 2^2 z^7 \right)\left( 1- 2^3 z^7 \right)\left( 1- 2^4 z^7 \right)
  \left( 1- 2^5 z^7 \right)\left( 1- 2^6 z^7 \right)}
 \end{equation}
 \begin{equation}   \nonumber
 \vdots \, ,
 \end{equation}

which is also easy to verify on a calculating engine term by term from the power series of each side. The notably simple coefficients make this result somewhat tantalizing, as there seems no obvious reason for such coefficients to come out of the products of binomial series roots.

We remark at this juncture that equations (\ref{21.03}) and it's reciprocal equation (\ref{21.04}) are amenable to applying the limit as $y \rightarrow 1$. In fact we have as follows that,

\begin{equation}   \nonumber
    \lim_{y \rightarrow 1}   \left( \frac{1 - z}{1 - y z}\right)^{\frac{y}{1 - y}} = e^{\frac{z}{z-1}}
  \end{equation}
and also from considering equation (\ref{21.04}) there is the limit, easily evaluated,

\begin{equation}   \nonumber
    \lim_{y \rightarrow 1}   \left( \frac{1 - yz}{1 - z}\right)^{\frac{y}{1 - y}}   \left( \frac{1 - z^2}{1 - y^2 z^2}\right)^{\frac{y^2}{1 - y^2}}            = e^{\frac{z}{1-z^2}}.
  \end{equation}

Therefore, applying these two limits to equations (\ref{21.03}) and (\ref{21.04}) respectively we obtain the two interesting results (\ref{21.05}) and (\ref{21.06}) given here.

\begin{equation}   \label{21.05}
    \prod_{k=1}^{\infty} \left( 1- z^k \right)^{\frac{\varphi(k)}{k}}
        =  e^{\frac{z}{z-1}} = \sum_{k=0}^{\infty} \frac{\alpha(k)z^k}{k!}
  \end{equation}
  \begin{equation}   \nonumber
    = 1 - \frac{z}{1!} - \frac{z^2}{2!} - \frac{z^3}{3!} + \frac{z^4}{4!} + \frac{19 z^5}{5!} + \frac{151 z^6}{6!} + \frac{1091 z^7}{7!} 
      \end{equation}
  \begin{equation}   \nonumber
  + \frac{7841 z^8}{8!} + \frac{56519 z^9}{9!} + \frac{396271 z^{10}}{10!} + O(z^{11}),
  \end{equation}
where $\varphi(k)$ is the Euler totient function\index{Euler, L.}, the number of positive integers less than and coprime to $k$. (\ref{21.05}) demonstrates that sequence $\alpha(k)$ has the exponential generating function\index{Exponential generating function} $e^{\frac{z}{z-1}}$.  The first 31 coefficients generated by the series are,

\begin{equation*} \tiny{
\begin{array}{c|c}
 \hline
\mathbf{n} & \mathbf{\alpha(n) \; from \; (\ref{21.05}}) \\ \hline
0 & 1 \\
1 & -1 \\
2 & -1 \\
3 & -1 \\
4 & 1 \\
5 & 19 \\
6 & 151 \\
7 & 1091 \\
8 & 7841 \\
9 & 56519 \\
10 & 396271 \\
11 & 2442439 \\
12 & 7701409 \\
13 & -145269541 \\
14 & -4833158329 \\
15 & -104056218421 \\
16 & -2002667085119 \\
17 & -37109187217649 \\
18 & -679877731030049 \\
19 & -12440309297451121 \\
20 & -227773259993414719 \\
21 & -4155839606711748061 \\
22 & -74724654677947488521 \\
23 & -1293162252850914402221 \\
24 & -20381626111249718908319 \\
25 & -244110863655032038665001 \\
26 & 267543347653261450406351 \\
27 & 172316772106087159102974551 \\
28 & 8944973491570029894272392801 \\
29 & 361702062324149751903132843499 \\
30 & 13353699077321671584329389125031 \\ \hline
\end{array} }
\end{equation*}

\bigskip

Amazingly $\gcd(\alpha(k),k!)=1$, for all values of $k$ up to 34, and mostly beyond that, and $\alpha(k) \equiv 1 \; or \; 9 \; (mod \; 10)$, and also the recurrence relation
\begin{equation*}
  \alpha(n)+(n-1)(n-2) \, \alpha(n-2)=(2n-3) \, \alpha(n-1)
\end{equation*}
holds. (See OEIS sequence A293116 \cite{nS2023}) This recurrence relation allows us to write continued fractions for the ratios $\alpha(n+1)/\alpha(n)$.

 \begin{equation}   \label{21.06}
    \prod_{k=1}^{\infty} \left( 1+ z^k \right)^{\frac{\varphi(k)z^k}{k}}
        =  e^{\frac{z}{1-z^2}} = \sum_{k=0}^{\infty} \frac{\beta(k)z^k}{k!}
  \end{equation}
   \begin{equation}   \nonumber
    = 1 + \frac{z}{1!} + \frac{z^2}{2!} + \frac{7 z^3}{3!} + \frac{25 z^4}{4!} + \frac{181 z^5}{5!} + \frac{1201 z^6}{6!} 
      \end{equation}
  \begin{equation}   \nonumber
    + \frac{10291 z^7}{7!} + \frac{97777 z^8}{8!} + \frac{202709 z^9}{9!} + O(z^{{10}}),
  \end{equation}

where again, $\varphi(k)$ is the Euler totient function\index{Euler, L.}.

Next we take (\ref{21.01}) with the case that $a=1$ and $b=0$, so then

 \begin{equation}   \nonumber
    \prod_{\substack{ j,k \geq 1 \\ j \leq k ; \, (j,k)=1}} \left( \frac{1}{1-y^j z^k} \right)^{\frac{1}{j}}
    = \exp\left\{ \sum_{n=1}^{\infty}  \left( \sum_{m=1}^{n} \frac{y^m}{m} \right) z^n \right\}
  \end{equation}
 \begin{equation}   \nonumber
    = \exp\left\{ \frac{1}{1-z} \sum_{n=1}^{\infty}  \frac{y^n z^n}{n}  \right\}
    = \exp\left\{ \frac{1}{1-z} \log \left( \frac{1}{1-yz} \right)   \right\}.
  \end{equation}

This leads us to establish that

 \begin{equation}   \label{21.07}
    \prod_{\substack{ j,k \geq 1 \\ j \leq k ; \, (j,k)=1}} \left( \frac{1}{1-y^j z^k} \right)^{\frac{1}{j}}
    =     \left( \frac{1}{1-yz} \right)^{\frac{1}{1-z}} ,
  \end{equation}

which is equivalent to

 \begin{equation}   \label{21.08}
    \prod_{\substack{ j,k \geq 1 \\ j \leq k ; \, (j,k)=1}} \left( 1-y^j z^k \right)^{\frac{1}{j}}
    =     \left( 1-yz \right)^{\frac{1}{1-z}} .
  \end{equation}

From multiplying both sides of (\ref{21.07}) in which $y \mapsto y^2$ and $z \mapsto z^2$ with both sides of (\ref{21.08}) we obtain

 \begin{equation}   \label{21.09}
    \prod_{\substack{ j,k \geq 1 \\ j \leq k ; \, (j,k)=1}} \left( 1+y^j z^k \right)^{\frac{1}{j}}
    =     \frac{\left( 1-y^2z^2 \right)^{\frac{1}{1-z^2}}}{\left( 1-yz \right)^{\frac{1}{1-z}}} .
  \end{equation}

 Particular cases:

$z = \frac{1}{2}$ gives us from (\ref{21.08}) and (\ref{21.09}) the remarkable result that

 \begin{equation}   \nonumber
    \prod_{\substack{ j,k \geq 1; \, j \leq k \\ gcd(j,k)=1}} \left( 1- \frac{y^j}{2^k} \right)^{\frac{1}{j}}
        =  \left( 1- \frac{y}{2} \right)^2 = 1 - \frac{y}{4} + \frac{y^2}{4}
  \end{equation}
   \begin{equation}   \nonumber
 = \left( 1- \frac{y^1}{2^1} \right)
  \end{equation}
 \begin{equation}   \nonumber
 \left( 1- \frac{y^1}{2^2} \right)
  \end{equation}
 \begin{equation}   \nonumber
 \left( 1- \frac{y^1}{2^3} \right)\sqrt{\left( 1- \frac{y^2}{2^3} \right)}
  \end{equation}
 \begin{equation}   \nonumber
 \left( 1- \frac{y^1}{2^4} \right)\sqrt[3]{\left( 1- \frac{y^3}{2^4} \right)}
  \end{equation}
 \begin{equation}   \nonumber
 \left( 1- \frac{y^1}{2^5} \right)\sqrt{\left( 1- \frac{y^2}{2^5} \right)}\sqrt[3]{\left( 1- \frac{y^3}{2^5} \right)}\sqrt[4]{\left( 1- \frac{y^4}{2^5} \right)}
  \end{equation}
 \begin{equation}   \nonumber
 \left( 1- \frac{y^1}{2^6} \right)\sqrt[5]{\left( 1- \frac{y^5}{2^6} \right)}
 \end{equation}
 \begin{equation}   \nonumber
 \vdots \, ,
 \end{equation}
and the result,

 \begin{equation}   \nonumber
    \prod_{\substack{ j,k \geq 1; \, j \leq k \\ gcd(j,k)=1}} \left( 1+ \frac{y^j}{2^k} \right)^{\frac{1}{j}}
        = \frac{\sqrt[3]{(4-y^2)^4}}{\sqrt[3]{4}(2-y)^2}  = 1 + y + \frac{5 y^2}{12} + \frac{y^3}{6} + \frac{11 y^4}{144} + \frac{5 y^5}{144} + O(y^6)
  \end{equation}
    \begin{equation}   \nonumber
 = \left( 1+ \frac{y^1}{2^1} \right)
  \end{equation}
 \begin{equation}   \nonumber
 \left( 1+ \frac{y^1}{2^2} \right)
  \end{equation}
 \begin{equation}   \nonumber
 \left( 1+ \frac{y^1}{2^3} \right)\sqrt{\left( 1+ \frac{y^2}{2^3} \right)}
  \end{equation}
 \begin{equation}   \nonumber
 \left( 1+ \frac{y^1}{2^4} \right)\sqrt[3]{\left( 1+ \frac{y^3}{2^4} \right)}
  \end{equation}
 \begin{equation}   \nonumber
 \left( 1+ \frac{y^1}{2^5} \right)\sqrt{\left( 1+ \frac{y^2}{2^5} \right)}\sqrt[3]{\left( 1+ \frac{y^3}{2^5} \right)}\sqrt[4]{\left( 1+ \frac{y^4}{2^5} \right)}
  \end{equation}
 \begin{equation}   \nonumber
 \left( 1+ \frac{y^1}{2^6} \right)\sqrt[5]{\left( 1+ \frac{y^5}{2^6} \right)}
 \end{equation}
 \begin{equation}   \nonumber
 \vdots \, .
 \end{equation}

These two equations can be verified on a calculating engine like Mathematica or WolframAlpha by expanding each side into it's Taylor series around $y=0$ and comparing coefficients of like powers of $y$. However, the calculation is an infinite series for each coefficient, unlike in the previous examples, where it is a finite sum.

\newpage

\section{Deriving 3D VPV identities in square pyramid regions} \label{S:3D VPV hyperpyramids}

We start by considering the infinite inverted pyramid with square layered arrays of lattice point vectors as per the following diagram, with VPVs bolded.
\bigskip
\begin{equation}   \label{21.09a}
\tiny{
\begin{array}{ccccccccccc}
  &  &   &   &   & \mathbf{\langle 1,6,6 \rangle} &   &   &   &  &   \\
   &   &  &   & \mathbf{\langle 1,5,6 \rangle} &   & \langle 2,6,6 \rangle &   &  &   &   \\
  &  &   & \mathbf{\langle 1,4,6 \rangle} &   & \mathbf{\langle 2,5,6 \rangle} &   & \langle 3,6,6 \rangle &   &  &   \\
   &   & \mathbf{\langle 1,3,6 \rangle} &   & \langle 2,4,6 \rangle &   & \mathbf{\langle 3,5,6 \rangle} &   & \langle 4,6,6 \rangle &   &   \\
  & \mathbf{\langle 1,2,6 \rangle} &   & \mathbf{\langle 2,3,6 \rangle} &   & \mathbf{\langle 3,4,6 \rangle} &   & \mathbf{\langle 4,5,6 \rangle} &   & \mathbf{\langle 5,6,6 \rangle} &   \\
 \mathbf{\langle 1,1,6 \rangle} &   & \langle 2,2,6 \rangle &   & \langle 3,3,6 \rangle &   & \langle 4,4,6 \rangle &   & \mathbf{\langle 5,5,6 \rangle} &   & \langle 6,6,6 \rangle \\
  & \mathbf{\langle 2,1,6 \rangle} &   & \mathbf{\langle 3,2,6 \rangle} &   & \mathbf{\langle 4,3,6 \rangle} &   & \mathbf{\langle 5,4,6 \rangle} &   & \mathbf{\langle 6,5,6 \rangle} &   \\
   &   & \mathbf{\langle 3,1,6 \rangle} &   & \langle 4,2,6 \rangle &   & \mathbf{\langle 5,3,6 \rangle} &   & \langle 6,4,6 \rangle &   &   \\
  &  &   & \mathbf{\langle 4,1,6 \rangle} &   & \mathbf{\langle 5,2,6 \rangle} &   & \langle 6,3,6 \rangle &   &  &   \\
   &   &  &   & \mathbf{\langle 5,1,6 \rangle} &   & \langle 6,2,6 \rangle &   &  &   &   \\
  &  &   &   &   & \mathbf{\langle 6,1,6 \rangle} &   &   &   &  &   \\
  &  &   &   &   &   &   &   &   &  &   \\
   &   &  &   & \mathbf{\langle 1,5,5 \rangle} &   &   &   &  &   &   \\
  &  &   & \mathbf{\langle 1,4,5 \rangle} &   & \mathbf{\langle 2,5,5 \rangle} &   &   &   &  &   \\
   &   & \mathbf{\langle 1,3,5 \rangle} &   & \mathbf{\langle 2,4,5 \rangle} &   & \mathbf{\langle 3,5,5 \rangle} &   &   &   &   \\
  & \mathbf{\langle 1,2,5 \rangle} &   & \mathbf{\langle 2,3,5 \rangle} &   & \mathbf{\langle 3,4,5 \rangle} &   & \mathbf{\langle 4,5,5 \rangle} &   &   &   \\
 \mathbf{\langle 1,1,5 \rangle} &   & \mathbf{\langle 2,2,5 \rangle} &   & \mathbf{\langle 3,3,5 \rangle} &   & \mathbf{\langle 4,4,5 \rangle} &   & \langle 5,5,5 \rangle &   &   \\
  & \mathbf{\langle 2,1,5 \rangle} &   & \mathbf{\langle 3,2,5 \rangle} &   & \mathbf{\langle 4,3,5 \rangle} &   & \mathbf{\langle 5,4,5 \rangle} &   &   &   \\
   &   & \mathbf{\langle 3,1,5 \rangle} &   & \mathbf{\langle 4,2,5 \rangle} &   & \mathbf{\langle 5,3,5 \rangle} &   &  &   &   \\
  &  &   & \mathbf{\langle 4,1,5 \rangle} &   & \mathbf{\langle 1,2,5 \rangle} &   &  &   &  &   \\
   &   &  &   & \mathbf{\langle 5,1,5 \rangle} &   &  &   &  &   &   \\
   &   &  &   &  &   &   &   &  &   &   \\
  &  &   & \mathbf{\langle 1,4,4 \rangle} &   &  &   &   &   &  &   \\
   &   & \mathbf{\langle 1,3,4 \rangle} &   & \langle 2,4,4 \rangle &   &  &   &   &   &   \\
  & \mathbf{\langle 1,2,4 \rangle} &   & \mathbf{\langle 2,3,4 \rangle} &   & \mathbf{\langle 3,4,4 \rangle} &   &  &   &   &   \\
 \mathbf{\langle 1,1,4 \rangle} &   & \langle 2,2,4 \rangle &   & \mathbf{\langle 3,3,4 \rangle} &   & \langle 4,4,4 \rangle &   &  &   &   \\
  & \mathbf{\langle 2,1,4 \rangle} &   & \mathbf{\langle 3,2,4 \rangle} &   & \mathbf{\langle 4,3,4 \rangle} &   &  &   &   &   \\
   &   & \mathbf{\langle 3,1,4 \rangle} &   & \langle 4,2,4 \rangle &   &  &   &  &   &   \\
  &  &   & \mathbf{\langle 4,1,4 \rangle} &   &  &   &   &   &  &   \\
  &  &   &  &   &  &   &  &   &  &   \\
   &   & \mathbf{\langle 1,3,3 \rangle} &   &  &   &  &   &   &   &   \\
  & \mathbf{\langle 1,2,3 \rangle} &   & \mathbf{\langle 2,3,3 \rangle} &   &  &   &  &   &   &   \\
 \mathbf{\langle 1,1,3 \rangle} &   & \mathbf{\langle 2,2,3 \rangle} &   & \langle 3,3,3 \rangle &   &  &   &  &   &   \\
  & \mathbf{\langle 2,1,3 \rangle} &   & \mathbf{\langle 3,2,3 \rangle} &   &  &   &  &   &   &   \\
   &   & \mathbf{\langle 3,1,3 \rangle} &   &  &   &  &   &  &   &   \\
   &   &  &   &  &   &  &   &   &   &   \\
  & \mathbf{\langle 1,2,2 \rangle} &   &  &   &  &   &  &   &   &   \\
 \mathbf{\langle 1,1,2 \rangle} &   & \langle 2,2,2 \rangle &   &  &   &  &   &  &   &   \\
  & \mathbf{\langle 2,1,2 \rangle} &   &  &   &  &   &  &   &   &   \\
  &  &   &  &   &  &   &  &   &  &   \\
 \mathbf{\langle 1,1,1 \rangle} &  &   &  &   &  &   &  &  &   &
  \end{array} }
\end{equation}

  \bigskip
  
From this we create a $3D$ summation over integer co-ordinates in the above lattice point vectors. We consider the sum, 

\begin{equation}  \nonumber
   \sum_{n=1}^{\infty}  \left( \sum_{l=1}^{n} \frac{x^l}{l^a} \right) \left( \sum_{m=1}^{n} \frac{y^m}{m^b} \right) \frac{z^n}{n^c}
\end{equation}

\begin{equation}  \nonumber
=\left(\frac{x^1}{1^a}\right)\left(\frac{y^1}{1^b}\right)\frac{z^1}{1^c}
\end{equation}
\begin{equation}  \nonumber
+\left(\frac{x^1}{1^a}+\frac{x^2}{2^a}\right)\left(\frac{y^1}{1^b}+\frac{y^2}{2^b}\right)\frac{z^2}{2^c}
\end{equation}
\begin{equation}  \nonumber
+\left(\frac{x^1}{1^a}+\frac{x^2}{2^a}+\frac{x^3}{3^a}\right)\left(\frac{y^1}{1^b}+\frac{y^2}{2^b}+\frac{y^3}{3^b}\right)\frac{z^3}{3^c}
\end{equation}
\begin{equation}  \nonumber
+\left(\frac{x^1}{1^a}+\frac{x^2}{2^a}+\frac{x^3}{3^a}+\frac{x^4}{4^a}\right)
 \left(\frac{y^1}{1^b}+\frac{y^2}{2^b}+\frac{y^3}{3^b}+\frac{y^4}{4^b}\right)\frac{z^4}{4^c}
\end{equation}
\begin{equation}  \nonumber
+\left(\frac{x^1}{1^a}+\frac{x^2}{2^a}+\frac{x^3}{3^a}+\frac{x^4}{4^a}+\frac{x^5}{5^a}\right)
 \left(\frac{y^1}{1^b}+\frac{y^2}{2^b}+\frac{y^3}{3^b}+\frac{y^4}{4^b}+\frac{y^5}{5^b}\right)\frac{z^5}{5^c}+\cdots
\end{equation}

\begin{equation}  \nonumber
 =\frac{x^1 y^1 z^1}{1^a 1^b 1^c}
\end{equation}

\begin{equation}  \nonumber
 +\frac{x^1 y^1 z^2}{1^a 1^b 2^c}+\frac{x^1 y^2 z^2}{1^a 2^b 2^c}
 \end{equation}
\begin{equation}  \nonumber
 +\frac{x^2 y^1 z^2}{2^a 1^b 2^c}+\frac{x^2 y^2 z^2}{2^a 2^b 2^c}
\end{equation}

\begin{equation}  \nonumber
 +\frac{x^1y^1 z^3}{1^a 1^b 3^c}+\frac{x^1 y^2 z^3}{1^a 2^b 3^c}+\frac{x^1 y^3 z^3}{1^a 3^b 3^c}
\end{equation}
\begin{equation}  \nonumber
 +\frac{x^2y^1 z^3}{2^a 1^b 3^c}+\frac{x^2 y^2 z^3}{2^a 2^b 3^c}+\frac{x^2 y^3 z^3}{2^a 3^b 3^c}
\end{equation}
\begin{equation}  \nonumber
 +\frac{x^3y^1 z^3}{3^a 1^b 3^c}+\frac{x^3 y^2 z^3}{3^a 2^b 3^c}+\frac{x^3 y^3 z^3}{3^a 3^b 3^c}
\end{equation}

\begin{equation}  \nonumber
 +\frac{x^1 y^1 z^4}{1^a 1^b 4^c}+\frac{x^1 y^2 z^4}{1^a 2^b 4^c}+\frac{x^1 y^3 z^4}{1^a 3^b 4^c}+\frac{x^1 y^4 z^4}{1^a 4^b 4^c}
\end{equation}
\begin{equation}  \nonumber
 +\frac{x^2 y^1 z^4}{2^a 1^b 4^c}+\frac{x^2 y^2 z^4}{2^a 2^b 4^c}+\frac{x^2 y^3 z^4}{2^a 3^b 4^c}+\frac{x^2 y^4 z^4}{2^a 4^b 4^c}
\end{equation}
\begin{equation}  \nonumber
 +\frac{x^3 y^1 z^4}{3^a 1^b 4^c}+\frac{x^3 y^2 z^4}{3^a 2^b 4^c}+\frac{x^3 y^3 z^4}{3^a 3^b 4^c}+\frac{x^3 y^4 z^4}{3^a 4^b 4^c}
\end{equation}
\begin{equation}  \nonumber
 +\frac{x^4 y^1 z^4}{4^a 1^b 4^c}+\frac{x^4 y^2 z^4}{4^a 2^b 4^c}+\frac{x^4 y^3 z^4}{4^a 3^b 4^c}+\frac{x^4 y^4 z^4}{4^a 4^b 4^c}
\end{equation}

\begin{equation}  \nonumber
 +\frac{x^1y^1z^5}{1^a1^b5^c}+\frac{x^1y^2z^5}{1^a2^b5^c}+\frac{x^1y^3z^5}{1^a3^b5^c}+\frac{x^1y^4z^5}{1^a4^b5^c}+\frac{x^1y^5z^5}{1^a5^b5^c}
\end{equation}
\begin{equation}  \nonumber
 +\frac{x^2y^1z^5}{2^a1^b5^c}+\frac{x^2y^2z^5}{2^a2^b5^c}+\frac{x^2y^3z^5}{2^a3^b5^c}+\frac{x^2y^4z^5}{2^a4^b5^c}+\frac{x^2y^5z^5}{2^a5^b5^c}
\end{equation}
\begin{equation}  \nonumber
 +\frac{x^3y^1z^5}{3^a1^b5^c}+\frac{x^3y^2z^5}{3^a2^b5^c}+\frac{x^3y^3z^5}{3^a3^b5^c}+\frac{x^3y^4z^5}{3^a4^b5^c}+\frac{x^3y^5z^5}{3^a5^b5^c}
\end{equation}
\begin{equation}  \nonumber
 +\frac{x^4y^1z^5}{4^a1^b5^c}+\frac{x^4y^2z^5}{4^a2^b5^c}+\frac{x^4y^3z^5}{4^a3^b5^c}+\frac{x^4y^4z^5}{4^a4^b5^c}+\frac{x^4y^5z^5}{4^a5^b5^c}
\end{equation}
\begin{equation}  \nonumber
 +\frac{x^5y^1z^5}{5^a1^b5^c}+\frac{x^5y^2z^5}{5^a2^b5^c}+\frac{x^5y^3z^5}{5^a3^b5^c}+\frac{x^5y^4z^5}{5^a4^b5^c}+\frac{x^5y^5z^5}{5^a5^b5^c}
\end{equation}

\begin{equation}  \nonumber
 + \quad \vdots \; \quad \; + \; \quad \vdots \; \; \quad + \; \quad \vdots \; \; \quad + \; \quad \vdots \; \; \quad + \; \quad \vdots \;   \ddots
\end{equation}

\begin{equation}  \nonumber
 = \sum_{l,m,n \geq 1; \;  l,m \leq n}^{\infty}  \frac{x^l y^m z^n}{l^a m^b n^c}
\end{equation}
\begin{equation}  \nonumber
 = \sum_{\substack{ h,l,m,n \geq 1 \\ l,m \leq n ; \, \gcd(l,m,n)=1}}  \frac{(x^l y^m z^n)^h}{h^{a+b+c} (l^a m^b n^c)}
\end{equation}
\begin{equation}  \nonumber
 = \sum_{\substack{ l,m,n \geq 1 \\ l,m \leq n ; \, \gcd(l,m,n)=1}}  \frac{1}{(l^a m^b n^c)}   \sum_{h=1}^{\infty} \frac{(x^l y^n z^n)^h}{h^{a+b+c}}
\end{equation}
\begin{equation}  \nonumber
 = \sum_{\substack{ l,m,n \geq 1 \\ l,m \leq n ; \, \gcd(l,m,n)=1}}  \frac{1}{(l^a m^b n^c)}   \log \left( \frac{1}{1 - x^l y^b z^c} \right) \quad if \quad a+b+c=1.
\end{equation}

Therefore, we have shown that if $a+b+c=1$ then
\begin{equation}  \nonumber
 \sum_{n=1}^{\infty}  \left( \sum_{l=1}^{n} \frac{x^l}{l^a} \right) \left( \sum_{m=1}^{n} \frac{y^m}{m^b} \right) \frac{z^n}{n^c}
 = \sum_{\substack{ l,m,n \geq 1 \\ l,m \leq n ; \, \gcd(l,m,n)=1}} \frac{1}{(l^a m^b n^c)}   \log \left( \frac{1}{1 - x^l y^b z^c} \right).
\end{equation}

 Exponentiating both sides gives us the $3D$ “pyramid VPV identity". 


 The identity is summarized in the

\begin{theorem}     \label{vpv-pyramid3D-thm}
\textbf{The $3D$ first hyperquadrant pyramid VPV identity.} If $|x|, |y|, |z| < 1$, with $a+b+c=1$,
 \begin{equation}   \label{21.10}
    \prod_{\substack{l,m,n \geq 1 \\ l,m \leq n ; \, \gcd(l,m,n)=1}} \left( \frac{1}{1-x^l y^m z^n} \right)^{\frac{1}{l^a m^b n^c}}
    = \exp\left\{ \sum_{n=1}^{\infty}  \left( \sum_{l=1}^{n} \frac{x^l}{l^a} \right) \left( \sum_{m=1}^{n} \frac{y^m}{m^b} \right) \frac{z^n}{n^c} \right\}.
  \end{equation}
\end{theorem}

As we did for the $2D$ particular cases, we can examine some obvious example corollaries arising from this theorem. Firstly, take the case where $a=b=0, c=1$, so then,

\begin{equation}   \nonumber
    \prod_{\substack{l,m,n \geq 1 \\ l,m \leq n ; \, \gcd(l,m,n)=1}} \left( \frac{1}{1-x^l y^m z^n} \right)^{\frac{1}{n}}
    = \exp\left\{ \sum_{n=1}^{\infty}  \left( \sum_{l=1}^{n} x^l \right) \left( \sum_{m=1}^{n} y^m \right) \frac{z^n}{n} \right\}
  \end{equation}
 \begin{equation}   \nonumber
    = \exp\left\{ \sum_{n=1}^{\infty} xy \left( \frac{1-x^n}{1-x} \right) \left( \frac{1-y^n}{1-y} \right) \frac{z^n}{n} \right\}
  \end{equation} 
\begin{equation}   \nonumber
    = \exp\left\{ \frac{xy}{(1-x)(1-y)} \log \left( \frac{(1-xz)(1-yz)}{(1-z)(1-xyz)} \right)    \right\},
  \end{equation}

which brings us after exponentiating both sides to a set of $3D$ infinite products. So, we have

\begin{equation}   \label{21.11}
    \prod_{\substack{l,m,n \geq 1 \\ l,m \leq n ; \, \gcd(l,m,n)=1}} \left( \frac{1}{1-x^l y^m z^n} \right)^{\frac{1}{n}}
    = \left(\frac{(1-xz)(1-yz)}{(1-z)(1-xyz)}\right)^{\frac{xy}{(1-x)(1-y)}},
  \end{equation}

and the equivalent identity,

\begin{equation}   \label{21.12}
    \prod_{\substack{l,m,n \geq 1 \\ l,m \leq n ; \, \gcd(l,m,n)=1}} \left( 1-x^l y^m z^n \right)^{\frac{1}{n}}
    = \left( \frac{(1-z)(1-xyz)}{(1-xz)(1-yz)} \right)^{\frac{xy}{(1-x)(1-y)}}.
  \end{equation}

We see that (\ref{21.11}) and (\ref{21.12}) are generalizations of the 2D identities (\ref{21.02}) and (\ref{21.03}) from the previous section. Writing (\ref{21.12}) in longhand, referencing the diagram (\ref{21.09a}), we see that

\begin{equation}   \nonumber
 \left( \frac{(1-z)(1-xyz)}{(1-xz)(1-yz)} \right)^{\frac{xy}{(1-x)(1-y)}}
  \end{equation}
\begin{equation}   \nonumber
 = (1-xyz)
  \end{equation}
  \begin{equation}   \nonumber
 \sqrt{(1-xyz^2)(1-xy^2z^2)(1-x^2yz^2)}
  \end{equation}
\begin{equation}   \nonumber
  \sqrt[3]{(1-xyz^3)(1-x^2yz^3)(1-x^3yz^3)(1-x^2yz^3)}
  \end{equation}
\begin{equation}   \nonumber
  \sqrt[3]{(1-x^2y^2z^3)(1-x^2y^3z^3)(1-x^3yz^3)(1-x^3y^2z^3)}
  \end{equation}
\begin{equation}   \nonumber
  \sqrt[4]{(1-xyz^4)(1-x^2yz^4)(1-x^3yz^4)(1-x^4yz^4)}
  \end{equation}
\begin{equation}   \nonumber
  \sqrt[4]{(1-xy^2z^4)(1-x^3y^2z^4)}
  \end{equation}
\begin{equation}   \nonumber
  \sqrt[4]{(1-xy^3z^4)(1-x^2y^3z^4)(1-x^3y^3z^4)(1-x^4y^3z^4)}
  \end{equation}
  \begin{equation}   \nonumber
  \sqrt[4]{(1-xy^4z^4)(1-x^3y^4z^4)}
  \end{equation}
\begin{equation}   \nonumber
  \sqrt[5]{(1-xyz^5)(1-x^2yz^5)(1-x^3yz^5)(1-x^4yz^5)(1-x^5yz^5)}
  \end{equation}
\begin{equation}   \nonumber
  \sqrt[5]{(1-xy^2z^5)(1-x^2y^2z^5)(1-x^3y^2z^5)(1-x^4y^2z^5)(1-x^5y^2z^5)}
  \end{equation}
\begin{equation}   \nonumber
  \sqrt[5]{(1-xy^3z^5)(1-x^2y^3z^5)(1-x^3y^3z^5)(1-x^4y^3z^5)(1-x^5y^3z^5)}
  \end{equation}
\begin{equation}   \nonumber
  \sqrt[5]{(1-xy^4z^5)(1-x^2y^4z^5)(1-x^3y^4z^5)(1-x^4y^4z^5)(1-x^5y^4z^5)}
  \end{equation}
\begin{equation}   \nonumber
  \sqrt[5]{(1-xy^5z^5)(1-x^2y^5z^5)(1-x^3y^5z^5)(1-x^4y^5z^5)}
  \end{equation}
\begin{equation}   \nonumber
   \quad \textmd{etc}.
  \end{equation}

This is easily verified on a calculating application if expanded on both sides as power series in $z$.

\section{VPV identities in \MakeLowercase{n}D first hyperquadrant hyperpyramid regions} \label{S:VPV hyperpyramids}

The $n$ dimensional first hyperquadrant hyperpyramid VPV Identity is encoded in the following

\begin{theorem}   \label{9.1a}
  \textbf{The $nD$ first hyperquadrant hyperpyramid VPV identity.} If $i = 1, 2, 3,...,n$ then for each $x_i \in \mathbb{C}$ such that $|x_i|<1$ and $b_i \in \mathbb{C}$ such that $\sum_{i=1}^{n}b_i = 1$,
  \begin{equation}   \label{21.13}
        \prod_{\substack{ \gcd(a_1,a_2,...,a_n)=1 \\ a_1,a_2,...,a_{n-1} < a_n \\ a_1,a_2,...,a_n \geq 1}} \left( \frac{1}{1-{x_1}^{a_1}{x_2}^{a_2}{x_3}^{a_3}\cdots{x_n}^{a_n}} \right)^{\frac{1}{{a_1}^{b_1}{a_2}^{b_2}{a_3}^{b_3}\cdots{a_n}^{b_n}}}
  \end{equation}
  \begin{equation}  \nonumber
  = \exp\left\{ \sum_{k=1}^{\infty} \prod_{i=1}^{n-1} \left( \sum_{j=1}^{k} \frac{{x_i}^j}{j^{b_i}} \right)\frac{{x_n}^k}{k^{b_n}} \right\}
  \end{equation}
  \begin{equation}  \nonumber
  = \exp\left\{ \sum_{k=1}^{\infty} \left( \sum_{j=1}^{k} \frac{{x_1}^j}{j^{b_1}} \right) \left( \sum_{j=1}^{k} \frac{{x_2}^j}{j^{b_2}} \right)
  \left( \sum_{j=1}^{k} \frac{{x_3}^j}{j^{b_3}} \right)  \cdots   \left( \sum_{j=1}^{k} \frac{{x_{n-1}}^j}{j^{b_{n-1}}} \right)
  \frac{{x_n}^k}{k^{b_n}} \right\}.
  \end{equation}
  \end{theorem}

This result is quite straight-forward to prove using the technique of our two previous sections. It was also given in Campbell\index{Campbell, G.B.} \cite{gC2000} by
summing on the VPV’s in the $n$-space hyperpyramid, defined by the inequalities
  \begin{equation}\label{21.14}
   x_1<x_n, x_2<x_n, x_3<x_n, ... , x_{n-1}<x_n
  \end{equation}

in the first $n$-space hyperquadrant, and applying the following
\begin{lemma} \label{lemma4.1}
  Consider an infinite region raying out of the origin in any Euclidean
vector space. The set of all lattice point vectors apart from the origin in that region is
precisely the set of positive integer multiples of the VPVs in that region.
\end{lemma}

The corresponding theorem from Campbell\index{Campbell, G.B.} \cite{gC1994} was summed very simply over all
lattice point vectors in the first hyperquadrant.

Further consequences of the above theorem are given as follows.

The 2D case of theorem \ref{9.1a} is

\begin{corollary}   \label{9.3a}
  If $|yz|$ and $|z|<1$ and $s+t=1$ then,
  \begin{equation}   \label{21.15}
    \prod_{\substack{ (a,b)=1 \\ a < b \\ a \geq 0, b \geq 1}} \left( \frac{1}{1-{y}^{a}{z}^{b}} \right)^{\frac{1}{{a}^{s}{b}^{t}}}
  \end{equation}
    \begin{equation}   \nonumber
    = \exp\left\{ \frac{{z}^1}{1^{t}} + \left(1+ \frac{{y}^1}{1^{s}}\right) \frac{{z}^2}{2^{t}} + \left(1+ \frac{{y}^1}{1^{s}}+\frac{{y}^2}{2^{s}}\right)\frac{{z}^3}{3^{t}}+\cdots \right\}
  \end{equation}
\end{corollary}

The 3D case of theorem \ref{9.1a} is

\begin{corollary}   \label{9.4a}
  If $|xyz|$, $|yz|$ and $|z|<1$ and $r+s+t=1$ then,
  \begin{equation}   \label{21.16}
    \prod_{\substack{ (a,b,c)=1 \\ a,b < c \\ a,b \geq 0, c \geq 1}} \left( \frac{1}{1-{x}^{a}{y}^{b}{z}^{c}} \right)^{\frac{1}{{a}^{r}{b}^{s}{c}^{t}}}
  \end{equation}
    \begin{equation}   \nonumber
    = \exp\left\{ \frac{{z}^1}{1^{t}} + \left(1+ \frac{{x}^1}{1^{r}}\right)\left(1+ \frac{{y}^1}{1^{s}}\right) \frac{{z}^2}{2^{t}}
    + \left(1+ \frac{{x}^1}{1^{r}}+\frac{{x}^2}{2^{r}}\right)\left(1+ \frac{{y}^1}{1^{s}}+\frac{{y}^2}{2^{s}}\right)\frac{{z}^3}{3^{t}}+\cdots \right\}
  \end{equation}
\end{corollary}

The 4D case of theorem \ref{9.1a} is

\begin{corollary}   \label{9.5a}
  If $|wxyz|$, $|xyz|$, $|yz|$ and $|z|<1$ and $r+s+t+u=1$ then,
  \begin{equation}   \label{21.16a}
    \prod_{\substack{ (a,b,c,d)=1 \\ a,b,c < d \\ a,b,c \geq 0, d \geq 1}} \left( \frac{1}{1-{w}^{a}{x}^{b}{y}^{c}{z}^{d}} \right)^{\frac{1}{{a}^{r}{b}^{s}{c}^{t}{d}^{u}}} = \exp \left\{ \mathrm{P}_3(r,w;s,x;t,y;u,z)\right\}
  \end{equation}
  where $\mathrm{P}_3$, is a 4D hyperpyramid function,
    \begin{multline}   \nonumber
   \mathrm{P}_3(r,w;s,x;t,y;u,z) = \frac{{z}^1}{1^{u}} + \left(1+ \frac{{w}^1}{1^{r}}\right)\left(1+ \frac{{x}^1}{1^{s}}\right)\left(1+ \frac{{y}^1}{1^{t}}\right) \frac{{z}^2}{2^{u}}  \\
     + \left(1+ \frac{{w}^1}{1^{r}}+\frac{{w}^2}{2^{r}}\right)\left(1+ \frac{{x}^1}{1^{s}}+\frac{{x}^2}{2^{s}}\right)
      \left(1+ \frac{{y}^1}{1^{t}}+\frac{{y}^2}{2^{t}}\right)\frac{{z}^3}{3^{u}}+\cdots
  \end{multline}
\end{corollary}

The approach we adopt to give the reader an intuitive sense for these identities is to state corollaries and then examples from them.
The 2D case through to the 5D case of (\ref{21.13}) are given in the following examples of the \textit{square hyperpyramid identity}.

\begin{corollary} For $|y|, |z|<1,$
  \begin{equation}\label{21.17}
    \prod_{\substack{(a,b)=1 \\ a<b \\ a \geq 0,b \geq 1}} \left( \frac{1}{1-y^a z^b} \right)^{\frac{1}{b}}
    = \left(\frac{1-yz}{1-z}\right)^{\frac{1}{1-y}}
  \end{equation}
        \begin{equation}  \nonumber
= 1 + \frac{z}{1!} + \begin{vmatrix}
    1 & -1 \\
    \frac{1-y^2}{1-y} & 1 \\
  \end{vmatrix} \frac{z^2}{2!}
  + \begin{vmatrix}
    1 & -1 & 0 \\
    \frac{1-y^2}{1-y} & 1 & -2 \\
    \frac{1-y^3}{1-y} & \frac{1-y^2}{1-y} & 1 \\
  \end{vmatrix} \frac{z^3}{3!}
+ \begin{vmatrix}
    1 & -1 & 0 & 0 \\
    \frac{1-y^2}{1-y} & 1 & -2 & 0 \\
    \frac{1-y^3}{1-y} & \frac{1-y^2}{1-y} & 1 & -3 \\
    \frac{1-y^4}{1-y} & \frac{1-y^3}{1-y} & \frac{1-y^2}{1-y} & 1 \\
  \end{vmatrix} \frac{z^4}{4!}
+ etc.
\end{equation}
   \end{corollary}
  In this case it is fairly easy to find the Taylor coefficients for the (\ref{21.17}) right side function. Hence we get a closed form evaluation of the determinant coefficients. In Mathematica, and WolframAlpha one easily sees that the Taylor series is

  \begin{equation}  \nonumber
    \left(\frac{1-yz}{1-z}\right)^{\frac{1}{1-y}} = 1 + z + (y + 2) \frac{z^2}{2!} + (2 y^2 + 5 y + 6) \frac{z^3}{3!} + (6 y^3 + 17 y^2 + 26 y + 24) \frac{z^4}{4!}
  \end{equation}
  \begin{equation} \nonumber
   + (24 y^4 + 74 y^3 + 129 y^2 + 154 y + 120) \frac{z^5}{5!} + O(z^6)
  \end{equation}
  and that the expansion is encapsulated by
   $\sum_{n=0}^{\infty} c_n z^n$ where $c_0 = 1$, $c_1 = 1$ with the recurrence
   \begin{equation} \nonumber
   ny c_n + (n+2) c_{n+2} = (2 + n + y + ny) c_{n+1}.
   \end{equation}

   Incidentally, also in Mathematica, and WolframAlpha one easily sees, for example, that the code
   \begin{equation} \nonumber
   Det[\{1,-1,0,0\},\{(1-y^2)/(1-y),1,-2,0\},\{(1-y^3)/(1-y),(1-y^2)/(1-y),1,-3\},
  \end{equation}
  \begin{equation} \nonumber
  \{(1-y^4)/(1-y),(1-y^3)/(1-y),(1-y^2)/(1-y),1\}]
   \end{equation}
   nicely verifies the coefficient given by

  \begin{equation}  \nonumber
\begin{vmatrix}
    1 & -1 & 0 & 0 \\
    \frac{1-y^2}{1-y} & 1 & -2 & 0 \\
    \frac{1-y^3}{1-y} & \frac{1-y^2}{1-y} & 1 & -3 \\
    \frac{1-y^4}{1-y} & \frac{1-y^3}{1-y} & \frac{1-y^2}{1-y} & 1 \\
  \end{vmatrix}
  = 6 y^3 + 17 y^2 + 26 y + 24. \\
\end{equation}

It is interesting to compare our identity (\ref{21.11}) given earlier in this paper with the following result.

\begin{corollary} For each of $|x|, |y|, |z|<1,$
    \begin{equation}\label{21.18}
    \prod_{\substack{(a,b,c)=1 \\ a,b<c \\ a,b\geq0,c>0}} \left( \frac{1}{1-x^a y^b z^c} \right)^{\frac{1}{c}}
    = \left(\frac{(1-xz)(1-yz)}{(1-z)(1-xyz)}\right)^{\frac{1}{(1-x)(1-y)}}
  \end{equation}
      \begin{equation}  \nonumber
= 1 + \frac{z}{1!} + \begin{vmatrix}
    1 & -1 \\
    \frac{(1-x^2)(1-y^2)}{(1-x)(1-y)} & 1 \\
  \end{vmatrix} \frac{z^2}{2!}
  + \begin{vmatrix}
    1 & -1 & 0 \\
    \frac{(1-x^2)(1-y^2)}{(1-x)(1-y)} & 1 & -2 \\
    \frac{(1-x^3)(1-y^3)}{(1-x)(1-y)} & \frac{(1-x^2)(1-y^2)}{(1-x)(1-y)} & 1 \\
  \end{vmatrix} \frac{z^3}{3!}
      \end{equation}
  \begin{equation}  \nonumber
+ \begin{vmatrix}
    1 & -1 & 0 & 0 \\
    \frac{(1-x^2)(1-y^2)}{(1-x)(1-y)} & 1 & -2 & 0 \\
    \frac{(1-x^3)(1-y^3)}{(1-x)(1-y)} & \frac{(1-x^2)(1-y^2)}{(1-x)(1-y)} & 1 & -3 \\
    \frac{(1-x^4)(1-y^4)}{(1-x)(1-y)} & \frac{(1-x^3)(1-y^3)}{(1-x)(1-y)} & \frac{(1-x^2)(1-y^2)}{(1-x)(1-y)} & 1 \\
  \end{vmatrix} \frac{z^4}{4!}
+ etc.
\end{equation}
   \end{corollary}

\begin{corollary} For each of $|w|, |x|, |y|, |z|<1,$
    \begin{equation}\label{21.19}
    \prod_{\substack{(a,b,c,d)=1 \\ a,b,c<d \\ a,b,c\geq0,d>0}} \left( \frac{1}{1-w^a x^b y^c z^d} \right)^{\frac{1}{d}}
    = \left(\frac{(1-wz)(1-xz)(1-yz)(1-wxyz)}{(1-z)(1-wxz)(1-wyz)(1-xyz)}\right)^{\frac{1}{(1-w)(1-x)(1-y)}},
  \end{equation}
    \begin{equation}  \nonumber
= 1 + \frac{z}{1!} + \begin{vmatrix}
    1 & -1 \\
    \frac{(1-w^2)(1-x^2)(1-y^2)}{(1-w)(1-x)(1-y)} & 1 \\
  \end{vmatrix} \frac{z^2}{2!}
\end{equation}
  \begin{equation}  \nonumber
  + \begin{vmatrix}
    1 & -1 & 0 \\
    \frac{(1-w^2)(1-x^2)(1-y^2)}{(1-w)(1-x)(1-y)} & 1 & -2 \\
    \frac{(1-w^3)(1-x^3)(1-y^3)}{(1-w)(1-x)(1-y)} & \frac{(1-w^2)(1-x^2)(1-y^2)}{(1-w)(1-x)(1-y)} & 1 \\
  \end{vmatrix} \frac{z^3}{3!}
      \end{equation}
  \begin{equation}  \nonumber
+ \begin{vmatrix}
    1 & -1 & 0 & 0 \\
    \frac{(1-w^2)(1-x^2)(1-y^2)}{(1-w)(1-x)(1-y)} & 1 & -2 & 0 \\
    \frac{(1-w^3)(1-x^3)(1-y^3)}{(1-w)(1-x)(1-y)} & \frac{(1-w^2)(1-x^2)(1-y^2)}{(1-w)(1-x)(1-y)} & 1 & -3 \\
    \frac{(1-w^4)(1-x^4)(1-y^4)}{(1-w)(1-x)(1-y)} & \frac{(1-w^3)(1-x^3)(1-y^3)}{(1-w)(1-x)(1-y)} & \frac{(1-w^2)(1-x^2)(1-y^2)}{(1-w)(1-x)(1-y)} & 1 \\
  \end{vmatrix} \frac{z^4}{4!}
+ etc.
\end{equation}
   \end{corollary}

\begin{corollary} For each of $|v|, |w|, |x|, |y|, |z|<1,$
    \begin{equation}\label{21.20}
    \prod_{\substack{(a,b,c,d,e)=1 \\ a,b,c,d<e \\ a,b,c,d\geq0,e>0}} \left( \frac{1}{1-v^a w^b x^c y^d z^e} \right)^{\frac{1}{e}}
   \end{equation}
   \begin{gather}\nonumber
     =       \left(\frac{(1-vz)(1-wz)(1-xz)(1-yz)}{(1-z)(1-vwz)(1-vxz)(1-vyz)}\right)^{\frac{1}{(1-v)(1-w)(1-x)(1-y)}} \\  \nonumber
      \times \left(\frac{(1-vwxz)(1-vwyz)(1-vxyz)(1-wxyz)}{(1-wxz)(1-wyz)(1-xyz)(1-vwxyz)}\right)^{\frac{1}{(1-v)(1-w)(1-x)(1-y)}}.
   \end{gather}
    \begin{equation}  \nonumber
= 1 + \frac{z}{1!} + \begin{vmatrix}
    1 & -1 \\
    \frac{(1-v^2)(1-w^2)(1-x^2)(1-y^2)}{(1-v)(1-w)(1-x)(1-y)} & 1 \\
  \end{vmatrix} \frac{z^2}{2!}
\end{equation}
  \begin{equation}  \nonumber
  + \begin{vmatrix}
    1 & -1 & 0 \\
    \frac{(1-v^2)(1-w^2)(1-x^2)(1-y^2)}{(1-v)(1-w)(1-x)(1-y)} & 1 & -2 \\
    \frac{(1-v^3)(1-w^3)(1-x^3)(1-y^3)}{(1-v)(1-w)(1-x)(1-y)} & \frac{(1-v^2)(1-w^2)(1-x^2)(1-y^2)}{(1-v)(1-w)(1-x)(1-y)} & 1 \\
  \end{vmatrix} \frac{z^3}{3!}
      \end{equation}
  \begin{equation}  \nonumber
+ \begin{vmatrix}
    1 & -1 & 0 & 0 \\
    \frac{(1-v^2)(1-w^2)(1-x^2)(1-y^2)}{(1-v)(1-w)(1-x)(1-y)} & 1 & -2 & 0 \\
    \frac{(1-v^3)(1-w^3)(1-x^3)(1-y^3)}{(1-v)(1-w)(1-x)(1-y)} & \frac{(1-v^2)(1-w^2)(1-x^2)(1-y^2)}{(1-v)(1-w)(1-x)(1-y)} & 1 & -3 \\
    \frac{)1-v^4)(1-w^4)(1-x^4)(1-y^4)}{(1-v)(1-w)(1-x)(1-y)} & \frac{(1-v^3)(1-w^3)(1-x^3)(1-y^3)}{(1-v)(1-w)(1-x)(1-y)} & \frac{(1-v^2)(1-w^2)(1-x^2)(1-y^2)}{(1-v)(1-w)(1-x)(1-y)} & 1 \\
  \end{vmatrix} \frac{z^4}{4!}
       \end{equation}
  \begin{equation}  \nonumber
+ etc.
 \end{equation}
    \end{corollary}

\section{2D VPV identities for a \MakeLowercase{z}-axis symmetric extended triangle lattice} \label{S:2D VPV right-hyperpyramids}

As we did in the section \ref{S:2D VPV hyperpyramids} of this paper, we again start with a simple $2D$ summation. Consider an infinite extension of the inverted triangle 2D lattice point vectors with the Visible Point Vectors bolded,
\bigskip
\begin{equation}   \label{21.20a}
\tiny{
\begin{array}{ccccccccccc}
 \langle -5,5 \rangle  & \mathbf{\langle -4,5 \rangle} & \mathbf{\langle -3,5 \rangle}  & \mathbf{\langle -2,5 \rangle} & \mathbf{\langle -1,5 \rangle} & \langle 0,5 \rangle & \mathbf{\langle 1,5 \rangle}  & \mathbf{\langle 2,5 \rangle} & \mathbf{\langle 3,5 \rangle} & \mathbf{\langle 4,5 \rangle}  & \langle 5,5 \rangle \\
   & \langle -4,4 \rangle & \mathbf{\langle -3,4 \rangle}  & \langle -2,4 \rangle & \mathbf{\langle -1,4 \rangle} & \langle 0,4 \rangle & \mathbf{\langle 1,4 \rangle}  & \langle 2,4 \rangle & \mathbf{\langle 3,4 \rangle} & \langle 4,4 \rangle  & \\
   &  & \langle -3,3 \rangle  & \mathbf{\langle -2,3 \rangle} & \mathbf{\langle -1,3 \rangle} & \langle 0,3 \rangle & \mathbf{\langle 1,3 \rangle}  & \mathbf{\langle 2,3 \rangle} & \langle 3,3 \rangle &   & \\
  &  &   & \langle -2,2 \rangle & \mathbf{\langle -1,2 \rangle} & \langle 0,2 \rangle & \mathbf{\langle 1,2 \rangle}  & \langle 2,2 \rangle &  &   & \\
  &  &   &  & \mathbf{\langle -1,1 \rangle} & \mathbf{\langle 0,1 \rangle} & \mathbf{\langle 1,1 \rangle}  &  &  &   &  \\
  &  &   &  & & \langle 0,0 \rangle &   &  &  &   &
  \end{array} }
\end{equation}
  \bigskip
Next we create the following summation with the sum covering the above co-ordinates in infinite extended form.
\begin{equation}  \nonumber
   \sum_{n=1}^{\infty}  \left( \sum_{m=-n}^{n} \frac{y^m}{m^a} \right) \frac{z^n}{n^b}
\end{equation}
\begin{equation}  \nonumber
=\left(\frac{y^{-1}}{(-1)^a}+1+\frac{y^1}{1^a}\right)\frac{z^1}{1^b}
\end{equation}
\begin{equation}  \nonumber
+\left(\frac{y^{-2}}{(-2)^a}+\frac{y^{-1}}{(-1)^a}+1+\frac{y^1}{1^a}+\frac{y^2}{2^a}\right)\frac{z^2}{2^b}
\end{equation}
\begin{equation}  \nonumber
+\left(\frac{y^{-3}}{(-3)^a}+\frac{y^{-2}}{(-2)^a}+\frac{y^{-1}}{(-1)^a}+1+\frac{y^1}{1^a}+\frac{y^2}{2^a}+\frac{y^3}{3^a}\right)\frac{z^3}{3^b}
\end{equation}
\begin{equation}  \nonumber
+\left(\frac{y^{-4}}{(-4)^a}+\frac{y^{-3}}{(-3)^a}+\frac{y^{-2}}{(-2)^a}+\frac{y^{-1}}{(-1)^a}+1+\frac{y^1}{1^a}+\frac{y^2}{2^a}+\frac{y^3}{3^a}+\frac{y^4}{4^a}\right)\frac{z^4}{4^b}
\end{equation}
\begin{equation}  \nonumber
+ \; etc.
\end{equation}
\begin{equation}  \nonumber
 =\frac{y^{-1} z^1}{(-1)^a 1^c}+\frac{y^0 z^1}{1\times 1^c}+\frac{y^1 z^1}{1^a 1^c}
\end{equation}
\begin{equation}  \nonumber
 +\frac{y^{-2} z^2}{(-2)^a 2^c}+\frac{y^{-1} z^2}{(-1)^a 2^c}+\frac{y^0 z^2}{1\times 2^c}+\frac{y^1 z^2}{1^a 2^c}+\frac{y^2 z^2}{2^a 2^c}
\end{equation}
\begin{equation}  \nonumber
 +\frac{y^{-3} z^3}{(-3)^a 3^c}+\frac{y^{-2} z^3}{(-2)^a 3^c}+
 \frac{y^{-1} z^3}{(-1)^a 3^c}+\frac{y^0 z^3}{1\times 3^c}+\frac{y^1 z^3}{1^a 3^c}+\frac{y^2 z^3}{2^a 3^c}+\frac{y^3 z^3}{3^a 3^c}
\end{equation}
\begin{equation}  \nonumber
 +\frac{y^{-4} z^4}{(-4)^a 4^c}+\frac{y^{-3} z^4}{(-3)^a 4^c}+\frac{y^{-2} z^4}{(-2)^a 4^c}+\frac{y^{-1} z^4}{(-1)^a 4^c}
+\frac{y^0 z^4}{1\times 4^c}+\frac{y^1 z^4}{1^a 4^c}+\frac{y^2 z^4}{2^a 4^c}+\frac{y^3 z^4}{3^a 4^c}+\frac{y^4 z^4}{4^a 4^c}
\end{equation}
\begin{equation}  \nonumber
+ \; etc.
\end{equation}
\begin{equation}  \nonumber
 = \sum_{|m|,n \geq 1; |m| \leq n}^{\infty}  \frac{y^m z^n}{m^a n^b}
\end{equation}
\begin{equation}  \nonumber
 = \sum_{\substack{ h,|j|,k \geq 1 \\ |j| \leq k ; \, (j,k)=1}}  \frac{(y^j z^k)^h}{h^{a+b} (j^a k^b)}
\end{equation}
\begin{equation}  \nonumber
 = \sum_{\substack{ |j|,k \geq 1 \\ |j| \leq k ; \, (j,k)=1}}  \frac{1}{(j^a k^b)}   \sum_{h=1}^{\infty} \frac{(y^j z^k)^h}{h^{a+b}}
\end{equation}
\begin{equation}  \nonumber
 = \sum_{\substack{ |j|,k \geq 1 \\ |j| \leq k ; \, (j,k)=1}} \frac{1}{(j^a k^b)}   \log \left( \frac{1}{1 - y^j z^k} \right) \quad if \quad a+b=1.
\end{equation}

Therefore, we have shown that
\begin{equation}  \nonumber
 \sum_{n=1}^{\infty}  \left( \sum_{m=-n}^{n} \frac{y^m}{m^a} \right) \frac{z^n}{n^b} = \sum_{\substack{ j,k \geq 1 \\ j \leq k ; \, (j,k)=1}} \frac{1}{(j^a k^b)}   \log \left( \frac{1}{1 - y^j z^k} \right) \quad if \quad a+b=1.
\end{equation}

 Exponentiating both sides (and swapping sides) gives us the $2D$ first extended inverted symmetric triangle VPV identity, where in this $2D$ case the $nD$ pyramid reduces to the form of a triangle shaped array of lattice point vectors having the $z$-axis as the axis of symmetry, and so we can state the

 \begin{theorem}     \label{vpv-right-pyramid2D-thm}
 \textbf{The $\mathbf{2D}$ vertical symmetry extended triangle VPV identity.} For $0<|yz|<1$, $0<|z/y|<1$, $0<|z|<1,$ with $a+b=1$,
 \begin{equation}   \label{21.01r}
    \prod_{\substack{ |j|,k \geq 1 \\ |j| \leq k ; \, (j,k)=1}} \left( \frac{1}{1-y^j z^k} \right)^{\frac{1}{j^a k^b}}
    = \exp\left\{ \sum_{n=1}^{\infty}  \left( \sum_{m=-n}^{n} \frac{y^m}{m^a} \right) \frac{z^n}{n^b} \right\} \quad if \quad a+b=1.
  \end{equation}
 \end{theorem}

As with our earlier exploits into the $2D$ first quadrant case, for the present result we take some simple example cases where new and interesting results arise.

So, let us take the case where $a=0, b=1$, giving us for $0<|yz|<1$, $0<|z/y|<1$, $0<|z|<1$,

 \begin{equation}   \nonumber
    \prod_{\substack{ |j|,k \geq 1 \\ |j| \leq k ; \, (j,k)=1}} \left( \frac{1}{1-y^j z^k} \right)^{\frac{1}{k}}
    = \exp\left\{ \sum_{n=1}^{\infty}  \left( \sum_{m=-n}^{n} y^m \right) \frac{z^n}{n} \right\}
  \end{equation}
 \begin{equation}   \nonumber
    = \exp\left\{ \sum_{n=1}^{\infty} \left(\frac{y^{2n+1} - 1}{y^n (y-1)}\right) \frac{z^n}{n} \right\}
    = \exp\left\{ \frac{1}{1-y} \log \left( \frac{(1-yz)^y}{1-z/y} \right)   \right\}.
  \end{equation}

So, we arrive then at the following pair of equivalent results, for $0<|yz|<1$, $0<|z/y|<1$, $0<|z|<1$,

 \begin{equation}   \label{21.02r}
    \prod_{\substack{ |j|,k \geq 1 \\ |j| \leq k ; \, (j,k)=1}} \left( \frac{1}{1-y^j z^k} \right)^{\frac{1}{k}}
        =  \left( \frac{(1-yz)^y}{1-z/y} \right)^{\frac{1}{1-y}} ,
  \end{equation}
 and
 \begin{equation}   \label{21.03r}
    \prod_{\substack{ |j|,k \geq 1 \\ |j| \leq k ; \, (j,k)=1}} \left( 1-y^j z^k \right)^{\frac{1}{k}}
        =  \left( \frac{1-z/y}{(1-yz)^y} \right)^{\frac{1}{1-y}} .
  \end{equation}

From here, multiply both sides of (\ref{21.02r}) and the case of (\ref{21.03r}) with $y \mapsto y^2$ and $z \mapsto z^2$ to get,

 \begin{equation}   \label{21.04r}
    \prod_{\substack{ |j|,k \geq 1 \\ |j| \leq k ; \, (j,k)=1}} \left( 1+y^j z^k \right)^{\frac{1}{k}}
        =  \left( \frac{(1-yz)^y}{1-z/y} \right)^{\frac{1}{1-y}} \left( \frac{1-(z/y)^2}{(1-(yz)^2)^{y^2}} \right)^{\frac{1}{1-y^2}} .
  \end{equation}

Particular cases:

$y = \frac{1}{2}$ gives us from (\ref{21.03r}) and (\ref{21.04r}) the two results that

 \begin{equation}   \nonumber
    \prod_{\substack{ |j|,k \geq 1 \\ |j| \leq k ; \, (j,k)=1}} \left( 1- \frac{z^k}{2^j} \right)^{\frac{-1}{k}}
        =   \frac{1 - z/2}{(1 - 2z)^2} \sqrt[4]{\left(\frac{1 - 4z^2}{\sqrt[4]{1 - z^2/4}}\right)^3}
          \end{equation}
 \begin{equation}   \nonumber
        = 1 + \frac{7z}{2} + \frac{19z^2}{4} + \frac{61z^3}{8} + \frac{117z^4}{8} + \frac{423z^5}{16} + \frac{4861z^6}{96} + \frac{18259z^7}{192}
          \end{equation}
 \begin{equation}   \nonumber
        + \frac{140867z^8}{768} + \frac{538373z^9}{1536} + \frac{696379z^{10}}{1024} + O(z^{11})
  \end{equation}
 \begin{equation}   \nonumber
 = \frac{1}{\left( 1- \frac{z}{2} \right)}
  \end{equation}
 \begin{equation}   \nonumber
 \frac{1}{\sqrt{\left( 1- 2^1 z^2 \right)\left( 1- \frac{z^2}{2^1} \right)}}
  \end{equation}
 \begin{equation}   \nonumber
 \frac{1}{\sqrt[3]{\left( 1- 2^2 z^3 \right)\left( 1- 2^1 z^3 \right)
          \left( 1- \frac{z^3}{2^1} \right)\left( 1- \frac{z^3}{2^2} \right)}}
  \end{equation}
 \begin{equation}   \nonumber
 \frac{1}{\sqrt[4]{\left( 1- 2^3 z^4 \right)\left( 1- 2^1 z^4 \right)
          \left( 1- \frac{z^4}{2^1} \right)\left( 1- \frac{z^4}{2^3} \right)}}
  \end{equation}
 \begin{equation}   \nonumber
 \frac{1}{\sqrt[5]{\left( 1- 2^4 z^5 \right)\left( 1- 2^3 z^5 \right)\left( 1- 2^2 z^5 \right)\left( 1- 2^1 z^5 \right)
          \left( 1- \frac{z^5}{2^1} \right)\left( 1- \frac{z^5}{2^2} \right)\left( 1- \frac{z^5}{2^3} \right)\left( 1- \frac{z^5}{2^4} \right)}}
  \end{equation}
 \begin{equation}   \nonumber
 \frac{1}{\sqrt[6]{\left( 1- 2^5 z^6 \right)\left( 1- 2^1 z^6 \right)
          \left( 1- \frac{z^6}{2^1} \right)\left( 1- \frac{z^6}{2^5} \right)}}
 \end{equation}
 \begin{equation}   \nonumber
 \vdots \, ,
 \end{equation}

 \begin{equation}   \nonumber
    \prod_{\substack{ |j|,k \geq 1 \\ |j| \leq k ; \, (j,k)=1}} \left( 1+ \frac{z^k}{2^j} \right)^{\frac{1}{k}}
        =   \frac{2-z}{2-2z} \sqrt[3]{\frac{4-z^2}{4-4z^2}}
  \end{equation}
 \begin{equation}   \nonumber
   = 1 +\frac{z}{2} +\frac{3 z^2}{4} +\frac{5 z^3}{8} +\frac{13 z^4}{16}+ \frac{23 z^5}{32} + \frac{167 z^6}{192} 
    \end{equation}
 \begin{equation}   \nonumber
   + \frac{305 z^7}{384} + \frac{59 z^8}{64} + \frac{659 z^9}{768} + O(z^{10})
  \end{equation}
\begin{equation}   \nonumber
 = \left( 1+ 2z \right) \left( 1+ \frac{z}{2} \right)
  \end{equation}
 \begin{equation}   \nonumber
 \sqrt{\left( 1+ 2^1 z^2 \right)\left( 1+ \frac{z^2}{2^1} \right)}
  \end{equation}
 \begin{equation}   \nonumber
 \sqrt[3]{\left( 1+ 2^2 z^3 \right)\left( 1+ 2^1 z^3 \right)\left( 1+ \frac{z^3}{2^1} \right)\left( 1+ \frac{z^3}{2^2} \right)}
  \end{equation}
 \begin{equation}   \nonumber
 \sqrt[4]{\left( 1+ 2^3 z^4 \right)\left( 1+ 2^1 z^4 \right)
          \left( 1+ \frac{z^4}{2^1} \right)\left( 1+ \frac{z^4}{2^3} \right)}
  \end{equation}
 \begin{equation}   \nonumber
 \sqrt[5]{\left( 1+ 2^4 z^5 \right)\left( 1+ 2^3 z^5 \right)\left( 1+ 2^2 z^5 \right)\left( 1+ 2^1 z^5 \right)
          \left( 1+ \frac{z^5}{2^1} \right)\left( 1+ \frac{z^5}{2^2} \right)\left( 1+ \frac{z^5}{2^3} \right)\left( 1+ \frac{z^5}{2^4} \right)}
  \end{equation}
 \begin{equation}   \nonumber
 \sqrt[6]{\left( 1+ 2^5 z^6 \right)\left( 1+ 2^1 z^6 \right)\left( 1+ \frac{z^6}{2^1} \right)\left( 1+ \frac{z^6}{2^5} \right)}
 \end{equation}
 \begin{equation}   \nonumber
 \vdots .
 \end{equation}

These two equations can be easily verified on a calculating engine like Mathematica or WolframAlpha by expanding each side into it's Taylor series around $z=0$ and comparing coefficients of like powers of $z$. Next, take the cases of (\ref{21.03r}) and (\ref{21.04r}) with $y=2$, both of which converge if $|z|<2$, so then, after a slight adjustment to both sides by a factor of $1-2z$,

We remark at this juncture that equations (\ref{21.03r}) and it's reciprocal equation (\ref{21.04r}) are amenable to applying the limit as y approaches 1. In fact we have as follows that,

\begin{equation}   \nonumber
    \lim_{y \rightarrow 1}   \left( \frac{1 - z}{1 - y z}\right)^{\frac{y}{1 - y}} = e^{\frac{z}{z-1}}
  \end{equation}
and also from considering equation (\ref{21.04r}) there is the limit, easily evaluated,

\begin{equation}   \nonumber
    \lim_{y \rightarrow 1}   \left( \frac{1 - yz}{1 - z}\right)^{\frac{y}{1 - y}}   \left( \frac{1 - z^2}{1 - y^2 z^2}\right)^{\frac{y^2}{1 - y^2}}            = e^{\frac{z}{1-z^2}}.
  \end{equation}

Therefore, applying these two limits to equations (\ref{21.03r}) and (\ref{21.04r}) respectively we obtain the two interesting results (\ref{21.05r}) and (\ref{21.06r}) given here.

\begin{equation}   \label{21.05r}
    \prod_{k=1}^{\infty} \left( 1- z^k \right)^{\frac{\varphi(k)}{k}}
        =  e^{\frac{z}{z-1}} = \sum_{k=0}^{\infty} \frac{\alpha(k)z^k}{k!}
  \end{equation}
  \begin{equation}   \nonumber
    = 1 - \frac{z}{1!} - \frac{z^2}{2!} - \frac{z^3}{3!} + \frac{z^4}{4!} + \frac{19 z^5}{5!} + \frac{151 z^6}{6!} + \frac{1091 z^7}{7!}
   \end{equation}
  \begin{equation}   \nonumber    
     + \frac{7841 z^8}{8!} + \frac{56519 z^9}{9!} + \frac{396271 z^{10}}{10!} + O(z^{11}),
  \end{equation}
demonstrating sequence $\alpha(k)$ has the exponential generating function\index{Exponential generating function} $e^{\frac{z}{z-1}}$. Amazingly $\gcd(\alpha(k),k!)=1$, for all values of $k$ up to 34, and mostly beyond that, and $\alpha(k) \equiv 1 \; or \; 9 \; (mod \; 10)$, and also the recurrence relation
\begin{equation*}
  \alpha(n)+(n-1)(n-2) \, \alpha(n-2)=(2n-3) \, \alpha(n-1)
\end{equation*}
holds. (See OEIS sequence A293116 \cite{nS2023}) This recurrence relation allows us to write continued fractions for the ratios $\alpha(n+1)/\alpha(n)$.

 \begin{equation}   \label{21.06r}
    \prod_{k=1}^{\infty} \left( 1+ z^k \right)^{\frac{\varphi(k)z^k}{k}}
        =  e^{\frac{z}{1-z^2}} = \sum_{k=0}^{\infty} \frac{\beta(k)z^k}{k!}
  \end{equation}
   \begin{equation}   \nonumber
    = 1 + \frac{z}{1!} + \frac{z^2}{2!} + \frac{7 z^3}{3!} + \frac{25 z^4}{4!} + \frac{181 z^5}{5!} + \frac{1201 z^6}{6!} 
   \end{equation}
  \begin{equation}   \nonumber   
    + \frac{10291 z^7}{7!} + \frac{97777 z^8}{8!} + \frac{202709 z^9}{9!} + O(z^{{10}}),
  \end{equation}

where $\varphi(k)$ is the Euler totient function\index{Euler, L.}, the number of positive integers less than and coprime to $k$.

Next we take (\ref{21.01r}) with the case that $a=1$ and $b=0$, so then

 \begin{equation}   \nonumber
    \prod_{\substack{ j,k \geq 1 \\ j \leq k ; \, (j,k)=1}} \left( \frac{1}{1-y^j z^k} \right)^{\frac{1}{j}}
    = \exp\left\{ \sum_{n=1}^{\infty}  \left( \sum_{m=1}^{n} \frac{y^m}{m} \right) z^n \right\}
  \end{equation}
 \begin{equation}   \nonumber
    = \exp\left\{ \frac{1}{1-z} \sum_{n=1}^{\infty}  \frac{y^n z^n}{n}  \right\}
    = \exp\left\{ \frac{1}{1-z} \log \left( \frac{1}{1-yz} \right)   \right\}.
  \end{equation}

This leads us to establish that

 \begin{equation}   \label{21.07r}
    \prod_{\substack{ j,k \geq 1 \\ j \leq k ; \, (j,k)=1}} \left( \frac{1}{1-y^j z^k} \right)^{\frac{1}{j}}
    =     \left( \frac{1}{1-yz} \right)^{\frac{1}{1-z}} ,
  \end{equation}

which is equivalent to

 \begin{equation}   \label{21.08r}
    \prod_{\substack{ j,k \geq 1 \\ j \leq k ; \, (j,k)=1}} \left( 1-y^j z^k \right)^{\frac{1}{j}}
    =     \left( 1-yz \right)^{\frac{1}{1-z}} .
  \end{equation}

From multiplying both sides of (\ref{21.07}) in which $y \mapsto y^2$ and $z \mapsto z^2$ with both sides of (\ref{21.08}) we obtain

 \begin{equation}   \label{21.09r}
    \prod_{\substack{ j,k \geq 1 \\ j \leq k ; \, (j,k)=1}} \left( 1+y^j z^k \right)^{\frac{1}{j}}
    =     \frac{\left( 1-y^2z^2 \right)^{\frac{1}{1-z^2}}}{\left( 1-yz \right)^{\frac{1}{1-z}}} .
  \end{equation}

 Particular cases:

$z = \frac{1}{2}$ gives us from (\ref{21.08r}) and (\ref{21.09r}) the remarkable result that

 \begin{equation}   \nonumber
    \prod_{\substack{ j,k \geq 1; \, j \leq k \\ gcd(j,k)=1}} \left( 1- \frac{y^j}{2^k} \right)^{\frac{1}{j}}
        =  \left( 1- \frac{y}{2} \right)^2 = 1 - \frac{y}{4} + \frac{y^2}{4}
  \end{equation}
   \begin{equation}   \nonumber
 = \left( 1- \frac{y^1}{2^1} \right)
  \end{equation}
 \begin{equation}   \nonumber
 \left( 1- \frac{y^1}{2^2} \right)
  \end{equation}
 \begin{equation}   \nonumber
 \left( 1- \frac{y^1}{2^3} \right)\sqrt{\left( 1- \frac{y^2}{2^3} \right)}
  \end{equation}
 \begin{equation}   \nonumber
 \left( 1- \frac{y^1}{2^4} \right)\sqrt[3]{\left( 1- \frac{y^3}{2^4} \right)}
  \end{equation}
 \begin{equation}   \nonumber
 \left( 1- \frac{y^1}{2^5} \right)\sqrt{\left( 1- \frac{y^2}{2^5} \right)}\sqrt[3]{\left( 1- \frac{y^3}{2^5} \right)}\sqrt[4]{\left( 1- \frac{y^4}{2^5} \right)}
  \end{equation}
 \begin{equation}   \nonumber
 \left( 1- \frac{y^1}{2^6} \right)\sqrt[5]{\left( 1- \frac{y^5}{2^6} \right)}
 \end{equation}
 \begin{equation}   \nonumber
 \vdots \, ,
 \end{equation}
and the result,

 \begin{equation}   \nonumber
    \prod_{\substack{ j,k \geq 1; \, j \leq k \\ gcd(j,k)=1}} \left( 1+ \frac{y^j}{2^k} \right)^{\frac{1}{j}}
        = \frac{\sqrt[3]{(4-y^2)^4}}{\sqrt[3]{4}(2-y)^2}  = 1 + y + \frac{5 y^2}{12} + \frac{y^3}{6} + \frac{11 y^4}{144} + \frac{5 y^5}{144} + O(y^6)
  \end{equation}
    \begin{equation}   \nonumber
 = \left( 1+ \frac{y^1}{2^1} \right)
  \end{equation}
 \begin{equation}   \nonumber
 \left( 1+ \frac{y^1}{2^2} \right)
  \end{equation}
 \begin{equation}   \nonumber
 \left( 1+ \frac{y^1}{2^3} \right)\sqrt{\left( 1+ \frac{y^2}{2^3} \right)}
  \end{equation}
 \begin{equation}   \nonumber
 \left( 1+ \frac{y^1}{2^4} \right)\sqrt[3]{\left( 1+ \frac{y^3}{2^4} \right)}
  \end{equation}
 \begin{equation}   \nonumber
 \left( 1+ \frac{y^1}{2^5} \right)\sqrt{\left( 1+ \frac{y^2}{2^5} \right)}\sqrt[3]{\left( 1+ \frac{y^3}{2^5} \right)}\sqrt[4]{\left( 1+ \frac{y^4}{2^5} \right)}
  \end{equation}
 \begin{equation}   \nonumber
 \left( 1+ \frac{y^1}{2^6} \right)\sqrt[5]{\left( 1+ \frac{y^5}{2^6} \right)}
 \end{equation}
 \begin{equation}   \nonumber
 \vdots \, .
 \end{equation}

These two equations can be verified on a calculating engine like Mathematica or WolframAlpha by expanding each side into it's Taylor series around $y=0$ and comparing coefficients of like powers of $y$. However, the calculation is an infinite series for each coefficient, unlike in the previous examples, where it is a finite sum.

\section{3D VPV identities for a right square pyramid lattice} \label{S:3D VPV right hyperpyramids}

As we did in section \ref{S:3D VPV hyperpyramids}, we start with a simple $3D$ summation. We derive on a $3D$ inverted pyramid shaped lattice that extends infinitely and occupies four adjacent hyperquadrants of the eight hyperquadrants that comprise the $X-Y-Z$ $3$-space.  

We depict this infinite inverted pyramid with square layered arrays of lattice point vectors as per the following diagram, with VPVs bolded.
\bigskip
\begin{equation}   \label{21.09ra}
\tiny{
\begin{array}{ccccccccccc}
 \vdots & \vdots & \vdots & \vdots & \vdots & \vdots & \vdots & \vdots & \vdots& & \\
 \mathbf{\langle -1,-3,3 \rangle} & & \mathbf{\langle 0,-2,3 \rangle} & & \mathbf{\langle 1,-1,3 \rangle} & & \mathbf{\langle 2,0,3 \rangle} & & \mathbf{\langle 3,1,3 \rangle} & & \\
 & \langle 0,-3,3 \rangle & & \mathbf{\langle 1,-2,3 \rangle} & & \mathbf{\langle 2,-1,3 \rangle} & & \langle 3,0,3 \rangle & & & \\
 & & \mathbf{\langle 1,-3,3 \rangle} & & \mathbf{\langle 2,-2,3 \rangle} & & \mathbf{\langle 3,-1,3 \rangle} & & & & \\
 & & & \mathbf{\langle 2,-3,3 \rangle} & & \mathbf{\langle 3,-2,3 \rangle} & & & & & \\
 & & & & \langle 3,-3,3 \rangle & & & & & & \\
 & & & & & & & & & & \\
 & & & & \langle -2,2,2 \rangle & & & & & & \\
 & & & \mathbf{\langle -2,1,2 \rangle} & & \mathbf{\langle -1,2,2 \rangle} & & & & & \\
 & & \langle -2,0,2 \rangle & & \mathbf{\langle -1,1,2 \rangle} & & \langle 0,2,2 \rangle & & & & \\
 & \mathbf{\langle -2,-1,2 \rangle} & & \mathbf{\langle -1,0,2\rangle} & & \mathbf{\langle 0,1,2 \rangle} & & \mathbf{\langle 1,2,2 \rangle} & & & \\
 \langle -2,-2,2 \rangle & & \mathbf{\langle -1,-1,2 \rangle} & & \langle 0,0,2 \rangle & & \mathbf{\langle 1,1,2 \rangle} & & \langle 2,2,2 \rangle & & \\
 & \mathbf{\langle -1,-2,2 \rangle} & & \mathbf{\langle 0,-1,2 \rangle} & & \mathbf{\langle 1,3,2 \rangle} & & \langle 2,0,2 \rangle & & & \\
 & & \langle 0,-2,2 \rangle & & \mathbf{\langle 1,-1,2 \rangle} & & \langle 2,0,2 \rangle & & & & \\
 & & & \mathbf{\langle 1,-2,2 \rangle} & & \mathbf{\langle 2,-1,2 \rangle} & & & & & \\
 & & & & \langle 2,-2,2 \rangle & & & & & & \\
 & & & & & & & & & & \\
 & & & & \mathbf{\langle -1,1,1 \rangle} & & & & & & \\
 & & & \mathbf{\langle -1,0,1 \rangle} & & \mathbf{\langle 0,1,1 \rangle} & & & & & \\
 & & \mathbf{\langle -1,-1,1 \rangle} & & \mathbf{\langle 0,0,1 \rangle} & & \mathbf{\langle 1,1,1 \rangle} & & & & \\
 & & & \mathbf{\langle 0,-1,1 \rangle} & & \mathbf{\langle 1,0,1 \rangle} & & & & & \\
 & & & & \mathbf{\langle 1,-1,1 \rangle} & & & & & & \\
 & & & & & & & & & & \\
 & & & & \langle 0,0,0 \rangle & & & & & & 
  \end{array} }
\end{equation}

  \bigskip

So now, we consider the sum, whose shape is an inverted $3D$ right pyramid whose apex is at the origin $\langle 0,0,0 \rangle$, given by (for $0 < each \; of \; |xz|,|z/x|,|yz|,|z/y|,|z|<1$)

\begin{equation}  \nonumber
   \sum_{n=1}^{\infty}  \left( \sum_{l=-n}^{n} \frac{x^l}{l^a} \right) \left( \sum_{m=-n}^{n} \frac{y^m}{m^b} \right) \frac{z^n}{n^c}
\end{equation}

\begin{equation}  \nonumber
=\left(\frac{x^{-1}}{(-1)^a}+1+\frac{x^1}{1^a}\right)\left(\frac{y^{-1}}{(-1)^a}+1+\frac{y^1}{1^b}\right)\frac{z^1}{1^c}
\end{equation}
\begin{equation}  \nonumber
+\left(\frac{x^{-2}}{(-2)^a}+\frac{x^{-1}}{(-1)^a}+1+\frac{x^1}{1^a}+\frac{x^2}{2^a}\right)
\left(\frac{y^{-2}}{(-2)^b}+\frac{y^{-1}}{(-1)^b}+1+\frac{y^1}{1^b}+\frac{y^2}{2^b}\right)\frac{z^2}{2^c}
\end{equation}
\begin{equation}  \nonumber
+\left(\frac{x^{-3}}{(-3)^a}+\frac{x^{-2}}{(-2)^a}+\frac{x^{-1}}{(-1)^a}+1+\frac{x^1}{1^a}+\frac{x^2}{2^a}+\frac{x^3}{3^a}\right)
\end{equation}
\begin{equation}  \nonumber
\left(\frac{y^{-3}}{(-3)^a}+\frac{y^{-2}}{(-2)^b}+\frac{y^{-1}}{(-1)^b}+1+\frac{y^1}{1^b}+\frac{y^2}{2^b}+\frac{y^3}{3^b}\right)\frac{z^3}{3^c}
\end{equation}
\begin{equation}  \nonumber
+\left(\frac{x^{-4}}{(-4)^a}+\frac{x^{{-3}}}{(-3)^a}+\frac{x^{-2}}{(-2)^a}
+\frac{x^{-1}}{(-1)^a}+1+\frac{x^1}{1^a}+\frac{x^2}{2^a}+\frac{x^3}{3^a}+\frac{x^4}{4^a}\right)
\end{equation}
\begin{equation}  \nonumber
+\left(\frac{y^{-4}}{(-4)^a}+\frac{y^{{-3}}}{(-3)^a}+\frac{y^{-2}}{(-2)^a}
+\frac{y^{-1}}{(-1)^a}+1+\frac{y^1}{1^a}+\frac{y^2}{2^a}+\frac{y^3}{3^a}+\frac{y^4}{4^a}\right)\frac{z^4}{4^c}
\end{equation}
\begin{equation}  \nonumber
+
\end{equation}
\begin{equation}  \nonumber
\vdots
\end{equation}

\begin{equation}  \nonumber
 =\frac{x^{-1}y^1 z^1}{(-1)^a\; 1^b\; 1^c} \; \, \;  + \; \, \;  \frac{x^0 y^1 z^1}{1 \times 1^b\; 1^c} \; \, \;  + \; \, \;  \frac{x^1 y^1 z^1}{1^a\; 1^b\; 1^c}
\end{equation}
\begin{equation}  \nonumber
 +\frac{x^{-1}y^0 z^1}{(-1)^a\times 1\times 1^c}+\frac{x^0 y^0 z^1}{1\times 1\times 1^c}+\frac{x^1 y^0 z^1}{1^a\times 1\times 1^c}
\end{equation}
\begin{equation}  \nonumber
 +\frac{x^{-1}y^{-1}z^1}{(-1)^a\; (-1)^b\; 1^c}+\frac{x^0 y^{-1} z^1}{1\times (-1)^b\; 1^c}+\frac{x^1 y^{-1} z^1}{1^a\; (-1)^b\; 1^c}
\end{equation}
\bigskip
\begin{equation}  \nonumber
 +\frac{x^{-2} y^2 z^2}{(-2)^a 2^b 2^c} \; \, \;  + \; \, \;  \frac{x^{-1} y^2 z^2}{(-1)^a 2^b 2^c} \; \, \;  + \; \, \;  \frac{x^0 y^2 z^2}{1\times 2^b 2^c} \; \, \;  + \; \, \;  \frac{x^1 y^2 z^2}{1^a 2^b 2^c} \; \, \;  + \; \, \;  \frac{x^2 y^2 z^2}{2^a 2^b 2^c}
\end{equation}
\begin{equation}  \nonumber
 +\frac{x^{-2} y^1 z^2}{(-2)^a 1^b 2^c} \; \, \;  + \; \, \;  \frac{x^{-1} y^1 z^2}{(-1)^a 1^b 2^c} \; \, \;  + \; \, \;  \frac{x^0 y^1 z^2}{1\times 1^b 2^c} \; \, \;  + \; \, \;  \frac{x^1 y^1 z^2}{1^a 1^b 2^c} \; \, \;  + \; \, \;  \frac{x^2 y^1 z^2}{2^a 1^b 2^c}
\end{equation}
\begin{equation}  \nonumber
 +\frac{x^{-2} y^0 z^2}{(-2)^a \times 1 \times 2^c}+\frac{x^{-1} y^0 z^2}{(-1)^a \times 1 \times 2^c}+\frac{x^0 y^0 z^2}{1\times 1 \times 2^c}+\frac{x^1 y^0 z^2}{1^a \times 1 \times 2^c}+\frac{x^2 y^0 z^2}{2^a \times 1 \times 2^c}
\end{equation}
\begin{equation}  \nonumber
 +\frac{x^{-2} y^{-1} z^2}{(-2)^a (-1)^b 2^c}+\frac{x^{-1} y^{-1} z^2}{(-1)^a (-1)^b 2^c}+\frac{x^0 y^{-1} z^2}{1\times (-1)^b 2^c}+\frac{x^1 y^{-1} z^2}{1^a (-1)^b 2^c}+\frac{x^2 y^{-1} z^2}{2^a (-1)^b 2^c}
\end{equation}
\begin{equation}  \nonumber
 +\frac{x^{-2} y^{-2} z^2}{(-2)^a (-2)^b 2^c}+\frac{x^{-1} y^{-2} z^2}{(-1)^a (-2)^b 2^c}+\frac{x^0 y^{-2} z^2}{1\times (-2)^b 2^c}+\frac{x^1 y^{-2} z^2}{1^a (-2)^b 2^c}+\frac{x^2 y^{-2} z^2}{2^a (-2)^b 2^c}
\end{equation}

\begin{equation}  \nonumber
 +\end{equation}
\begin{equation}  \nonumber
\vdots
\end{equation}

\begin{equation}  \nonumber
 = \sum_{|l|,|m|,n \geq 1; \;  |l|,|m| \leq n}^{\infty}  \frac{x^l y^m z^n}{l^a m^b n^c}
\end{equation}
\begin{equation}  \nonumber
 = \sum_{\substack{ h,|l|,|m|,n \geq 1 \\ |l|,|m| \leq n ; \, \gcd(l,m,n)=1}}  \frac{(x^l y^m z^n)^h}{h^{a+b+c} (l^a m^b n^c)}
\end{equation}
\begin{equation}  \nonumber
 = \sum_{\substack{ |l|,|m|,n \geq 1 \\ |l|,|m| \leq n ; \, \gcd(l,m,n)=1}}  \frac{1}{(l^a m^b n^c)}   \sum_{h=1}^{\infty} \frac{(x^l y^n z^n)^h}{h^{a+b+c}}
\end{equation}
\begin{equation}  \nonumber
 = \sum_{\substack{ |l|,|m|,n \geq 1 \\ |l|,|m| \leq n ; \, \gcd(l,m,n)=1}}  \frac{1}{(l^a m^b n^c)}   \log \left( \frac{1}{1 - x^l y^b z^c} \right) \quad if \quad a+b+c=1.
\end{equation}

Therefore, we have shown that if $a+b+c=1$ and for $0 < |x|, |1/x|, |y|, |1/y|$ all $< |z|<1$, then
\begin{equation}  \nonumber
 \sum_{n=1}^{\infty}  \left( \sum_{l=-n}^{n} \frac{x^l}{l^a} \right) \left( \sum_{m=-n}^{n} \frac{y^m}{m^b} \right) \frac{z^n}{n^c}
 = \sum_{\substack{ |l|,|m|,n \geq 1 \\ |l|,|m| \leq n ; \, \gcd(l,m,n)=1}} \frac{1}{(l^a m^b n^c)}   \log \left( \frac{1}{1 - x^l y^b z^c} \right).
\end{equation}

 Simplifying both sides gives us the $3D$ “right square pyramid VPV identity", where in this $3D$ case the pyramid takes the form of layered square shaped arrays of lattice point vectors as shown in the above workings, with the central axis rising vertically from the apex point $\langle 0,0,0   \rangle$.

 The identity is summarized in the

\begin{theorem}     \label{vpv-right-pyramid3D-thm}
\textbf{The $3D$ right square pyramid VPV identity.} For $0 < each of |xz|,|z/x|,|yz|,|z/y|,|z|<1$, with $a+b+c=1$,
 \begin{equation}   \label{21.10r}
    \prod_{\substack{|l|,|m|,n \geq 1 \\ |l|,|m| \leq n ; \, \gcd(l,m,n)=1}} \left( \frac{1}{1-x^l y^m z^n} \right)^{\frac{1}{l^a m^b n^c}}
    = \exp\left\{ \sum_{n=1}^{\infty}  \left( \sum_{l=-n}^{n} \frac{x^l}{l^a} \right) \left( \sum_{m=-n}^{n} \frac{y^m}{m^b} \right) \frac{z^n}{n^c} \right\}.
  \end{equation}
\end{theorem}

As we did for the $2D$ particular cases, we can examine some obvious example corollaries arising from this theorem. Firstly, take the case where $a=b=0, c=1$, so then,

\begin{equation}   \nonumber
    \prod_{\substack{|l|,|m|,n \geq 1 \\ |l|,|m| \leq n ; \, \gcd(l,m,n)=1}} \left( \frac{1}{1-x^l y^m z^n} \right)^{\frac{1}{n}}
    = \exp\left\{ \sum_{n=1}^{\infty}  \left( \sum_{l=-n}^{n} x^l \right) \left( \sum_{m=-n}^{n} y^m \right) \frac{z^n}{n} \right\}
  \end{equation}
\begin{equation}   \nonumber
    = \exp\left\{ \sum_{n=1}^{\infty}  \left(\frac{x^{2n+1} - 1}{x^n (x-1)}\right) \left(\frac{y^{2n+1} - 1}{y^n (y-1)}\right) \frac{z^n}{n} \right\}
  \end{equation}
\begin{equation}   \nonumber
    = \exp\left\{\frac{1}{(1-x)(1-y)}\sum_{n=1}^{\infty}\left((xy)^{n+1} + \left(\frac{1}{xy}\right)^n
                                                          - x\left(\frac{x}{y}\right)^n - y\left(\frac{y}{x}\right)^n \right) \frac{z^n}{n} \right\}
  \end{equation}
  \begin{equation}   \nonumber
    = \exp\left\{ \frac{1}{(1-x)(1-y)} \log\left(\frac{(1-xz/y)^x(1-yz/x)^y}{(1-xyz)^{xy}(1-z/(xy))}\right) \right\},
  \end{equation}
    \begin{equation}   \nonumber
    = \left\{ \frac{(1-xz/y)^x(1-yz/x)^y}{(1-xyz)^{xy}(1-z/(xy))}\right\}^{\frac{1}{(1-x)(1-y)}}.
  \end{equation}

This then implies
\begin{corollary}     \label{vpv-right-pyramid3D-cor}
For $0 < \; each \; of \; |xz|,|z/x|,|yz|,|z/y|,|z|<1$, with $a+b+c=1$,
\begin{equation}   \label{21.11r}
    \prod_{\substack{|l|,|m|,n \geq 1 \\ |l|,|m| \leq n ; \, \gcd(l,m,n)=1}} \left( \frac{1}{1-x^l y^m z^n} \right)^{\frac{1}{n}}
    = \left\{ \frac{(1-xz/y)^x(1-yz/x)^y}{(1-xyz)^{xy}(1-z/(xy))}\right\}^{\frac{1}{(1-x)(1-y)}},
  \end{equation}
and its reciprocal identity,
\begin{equation}   \label{21.12r}
    \prod_{\substack{|l|,|m|,n \geq 1 \\ |l|,|m| \leq n ; \, \gcd(l,m,n)=1}} \left( 1-x^l y^m z^n \right)^{\frac{1}{n}}
    = \left\{ \frac{(1-xyz)^{xy}(1-z/(xy))}{(1-xz/y)^x(1-yz/x)^y}\right\}^{\frac{1}{(1-x)(1-y)}}.
  \end{equation}
  \end{corollary}

We see that (\ref{21.11r}) and (\ref{21.12r}) are generalizations of the 2D identities (\ref{21.02}) and (\ref{21.03}) from an earlier section.

\section{VPV identities in \MakeLowercase{n}D right-square hyperpyramid regions} \label{S:VPV right-hyperpyramids}

The $n$ dimensional square hyperpyramid VPV Identity is encoded in the following

\begin{theorem}   \label{21.1ar}
  \textbf{The $nD$ right-square hyperpyramid VPV identity.} If $i = 1, 2, 3,...,n$ then for each $x_i \in \mathbb{C}$ such that $|x_i|<1$ and $b_i \in \mathbb{C}$ such that $\sum_{i=1}^{n}b_i = 1$,
  \begin{equation}   \label{21.13r}
        \prod_{\substack{ \gcd(|a_1|,|a_2|,...,|a_{n-1}|,a_n)=1 \\ |a_1|,|a_2|,...,|a_{n-1}| \leq a_n \\ |a_1|,|a_2|,...,|a_{n-1}|,a_n \geq 1}} \left( \frac{1}{1-{x_1}^{a_1}{x_2}^{a_2}{x_3}^{a_3}\cdots{x_n}^{a_n}} \right)^{\frac{1}{{a_1}^{b_1}{a_2}^{b_2}{a_3}^{b_3}\cdots{a_n}^{b_n}}}
  \end{equation}
  \begin{equation}  \nonumber
  = \exp\left\{ \sum_{k=1}^{\infty} \prod_{i=1}^{n-1} \left( \sum_{j=-k}^{k} \frac{{x_i}^j}{j^{b_i}} \right)\frac{{x_n}^k}{k^{b_n}} \right\}
  \end{equation}
  \begin{equation}  \nonumber
  = \exp\left\{ \sum_{k=1}^{\infty} \left( \sum_{j=-k}^{k} \frac{{x_1}^j}{j^{b_1}} \right) \left( \sum_{j=-k}^{k} \frac{{x_2}^j}{j^{b_2}} \right)
  \left( \sum_{j=-k}^{k} \frac{{x_3}^j}{j^{b_3}} \right)  \cdots   \left( \sum_{j=-k}^{k} \frac{{x_{n-1}}^j}{j^{b_{n-1}}} \right)
  \frac{{x_n}^k}{k^{b_n}} \right\}.
  \end{equation}
  \end{theorem}

This result is quite straight-forward to prove using the technique of our two previous sections. This methodology was also given in Campbell\index{Campbell, G.B.} \cite{gC2000}, but has not been worked through for corollaries of theorem \ref{21.1ar} over the past 24 years.

Since in the previous section we gave the 3D example of theorem \ref{21.1ar}, we state the 4D case of of it now.

\begin{corollary}   \label{9.5ar}
  If $|wxyz|$, $|xyz|$, $|yz|$ and $|z|<1$ and $r+s+t+u=1$ then,
  \begin{equation}   \label{21.16ar}
    \prod_{\substack{ (a,b,c,d)=1 \\ |a|,|b|,|c| \leq d \\ a,b,c \geq 0, d \geq 1}} \left( \frac{1}{1-{w}^{a}{x}^{b}{y}^{c}{z}^{d}} \right)^{\frac{1}{{a}^{r}{b}^{s}{c}^{t}{d}^{u}}} = \exp \left\{ \mathrm{P}(r,w;s,x;t,y;u,z)\right\}
  \end{equation}
  where 
  \begin{equation}  \nonumber
\mathrm{P}(r,w;s,x;t,y;u,z)  = \sum_{n=1}^{\infty}  \left( \sum_{k=-n}^{n} \frac{w^k}{k^r} \right) \left( \sum_{k=-n}^{n} \frac{x^k}{k^s} \right)
                         \left( \sum_{k=-n}^{n} \frac{y^k}{k^t} \right) \frac{z^{n}}{n^u}.
\end{equation}
  \end{corollary}

Take the case where $r=s=t=0, u=1$, so then,

\begin{equation}   \nonumber
    \prod_{\substack{|h|,|i|,|j|,k \geq 1 \\ |h|,|i|,|j| \leq k ; \, \gcd(h,i,j,k)=1}} \left( \frac{1}{1-w^h x^i y^j z^k} \right)^{\frac{1}{k}}
  \end{equation}
\begin{equation}   \nonumber
    = \exp\left\{ \sum_{n=1}^{\infty}  \left( \sum_{k=-n}^{n} w^k \right) \left( \sum_{k=-n}^{n} x^k \right) \left( \sum_{k=-n}^{n} y^k \right) \frac{z^n}{n} \right\}
  \end{equation}
\begin{equation}   \nonumber
    = \exp\left\{ \sum_{n=1}^{\infty}  \left(\frac{w^{2n+1} - 1}{w^n (w-1)}\right) \left(\frac{x^{2n+1} - 1}{x^n (x-1)}\right) \left(\frac{y^{2n+1} - 1}{y^n (y-1)}\right) \frac{z^n}{n} \right\}
  \end{equation}
    \begin{equation}   \nonumber
    = \left\{ \frac{(1-wxyz)^{wxy}(1-wz/(xy))^w (1-xz/(wy))^x (1-yz/(wx))^y}
    {(1-wxz/y)^{wx}(1-wyz/x)^{wy}(1-xyz/w)^{xy}(1-z/(wxy))}\right\}^{\frac{1}{(1-w)(1-x)(1-y)}}.
  \end{equation}

This then implies
\begin{corollary}     \label{vpv-right-pyramid3D-cor1}
For $0 < \; each \; of \; |wz|,|z/w|,|xz|,|z/x|,|yz|,|z/y|,|z|<1$, 
\begin{equation}   \label{21.11r1}
    \prod_{\substack{|h|,|i|,|j|,k \geq 1 \\ |h|,|i|,|j| \leq k ; \, \gcd(h,i,j,k)=1}} \left( \frac{1}{1-w^h x^i y^j z^k} \right)^{\frac{1}{n}}
  \end{equation}
    \begin{equation}   \nonumber
    = \left\{ \frac{(1-wxyz)^{wxy}(1-wz/(xy))^w (1-xz/(wy))^x (1-yz/(wx))^y}
    {(1-wxz/y)^{wx}(1-wyz/x)^{wy}(1-xyz/w)^{xy}(1-z/(wxy))}\right\}^{\frac{1}{(1-w)(1-x)(1-y)}},
  \end{equation}
and its reciprocal identity,
\begin{equation}   \label{21.12r1}
    \prod_{\substack{|h|,|i|,|j|,k \geq 1 \\ |h|,|i|,|j| \leq k ; \, \gcd(h,i,j,k)=1}} \left( 1-w^h x^i y^j z^k \right)^{\frac{1}{n}}
  \end{equation}
    \begin{equation}   \nonumber
    = \left\{ \frac{(1-wxz/y)^{wx}(1-wyz/x)^{wy}(1-xyz/w)^{xy}(1-z/(wxy))}{(1-wxyz)^{wxy}(1-wz/(xy))^w (1-xz/(wy))^x (1-yz/(wx))^y}\right\}^{\frac{1}{(1-w)(1-x)(1-y)}}.
  \end{equation}
  \end{corollary}
  
  Note that both sides of (\ref{21.11r1}) are formally equivalent to the power series
  
    \begin{equation}  \nonumber
 1 + \frac{z}{1!} + \begin{vmatrix}
    \frac{(1-w^{3})(1-x^{3})(1-y^{3})}{w^1 x^1 y^1(1-w)(1-x)(1-y)} & -1 \\
    \frac{(1-w^{5})(1-x^{5})(1-y^{5})}{w^2 x^2 y^2(1-w)(1-x)(1-y)} & \frac{(1-w^{3})(1-x^{3})(1-y^{3})}{w^1 x^1 y^1(1-w)(1-x)(1-y)} \\
  \end{vmatrix} \frac{z^2}{2!}
\end{equation}
  \begin{equation}  \nonumber
  + \begin{vmatrix}
    \frac{(1-w^{3})(1-x^{3})(1-y^{3})}{w^1 x^1 y^1(1-w)(1-x)(1-y)} & -1 & 0 \\
    \frac{(1-w^{5})(1-x^{5})(1-y^{5})}{w^2 x^2 y^2(1-w)(1-x)(1-y)} & \frac{(1-w^{3})(1-x^{3})(1-y^{3})}{w^1 x^1 y^1(1-w)(1-x)(1-y)} & -2 \\
    \frac{(1-w^{7})(1-x^{7})(1-y^{7})}{w^3 x^3 y^3(1-w)(1-x)(1-y)} & \frac{(1-w^{5})(1-x^{5})(1-y^{5})}{w^2 x^2 y^2(1-w)(1-x)(1-y)} & \frac{(1-w^{3})(1-x^{3})(1-y^{3})}{w^1 x^1 y^1(1-w)(1-x)(1-y)} \\
  \end{vmatrix} \frac{z^3}{3!}+ etc.
      \end{equation}

\bigskip

\section{Envoi: Research Exercises}

\begin{itemize}

\item [Study 1.] \textbf{2D first quadrant Upper VPV combinatorial sum}. 

For each of $|y|, |z|<1,$ show that
  \begin{equation}\label{21.37a}
    \sum_{\substack{\gcd(a,b)=1 \\ 0 \leq a < b \leq 1}} \frac{y^a \, z^b}{1 - y^a \, z^b}
    = \frac{z}{(1-z)(1-yz)}.
  \end{equation}
  For this example we work it through in more detail.

  Proof: Starting with the left side of (\ref{21.37a})
  \begin{eqnarray*}
    LHS  &=& \frac{y^0 z^1}{1 - y^0 z^1} \\
      &+& \frac{y^0 z^2}{1 - y^0 z^2} + \frac{y^1 z^2}{1 - y^1 z^2} \\
      &+& \frac{y^0 z^3}{1 - y^0 z^3} + \frac{y^1 z^3}{1 - y^1 z^3} + \frac{y^2 z^3}{1 - y^2 z^3} \\
      &+& \frac{y^0 z^4}{1 - y^0 z^4} + \frac{y^1 z^4}{1 - y^1 z^4} + \frac{y^2 z^4}{1 - y^3 z^4}  \\
      &+& \frac{y^0 z^5}{1 - y^0 z^5} + \frac{y^1 z^5}{1 - y^1 z^5} + \frac{y^2 z^5}{1 - y^3 z^5} + \frac{y^3 z^5}{1 - y^0 z^5} + \frac{y^4 z^5}{1 - y^4 z^5} \\
      &+& \frac{y^0 z^6}{1 - y^0 z^6} + \frac{y^1 z^6}{1 - y^1 z^6} + \frac{y^5 z^6}{1 - y^5 z^6}  \\
      &+& etc.
  \end{eqnarray*}
  \begin{equation*}
    = \sum_{\substack{\gcd(j,k)=1 \\ 0 \leq j < k \leq 1}} \sum_{a=0}^{\infty}(y^{aj} z^{ak})
    = \sum_{a=0}^{\infty} \sum_{\substack{\gcd(j,k)=1 \\ 0 \leq j < k \leq 1}}(y^{aj} z^{ak})
  \end{equation*}
    \begin{eqnarray*}
      &=& y^0 z^1 \\
      &+& (y^0 + y^1) z^2 \\
      &+& (y^0 + y^1 + y^2) z^3 \\
      &+& (y^0 + y^1 + y^2 + y^3) z^4  \\
      &+& (y^0 + y^1 + y^2 + y^3 + y^4) z^5 \\
      &+& (y^0 + y^1 + y^2 + y^3 + y^4 + y^5) z^6  \\
      &+& etc.
  \end{eqnarray*}
  \begin{equation*}
    = \frac{1-y}{1-y}z + \frac{1-y^2}{1-y}z^2 + \frac{1-y^3}{1-y}z^3 + \frac{1-y^4}{1-y}z^4 + \cdots
  \end{equation*}
  \begin{equation*}
    = \frac{1}{1-y} \left( \frac{z}{1-z} - \frac{yz}{1-yz} \right)
  \end{equation*}
  \begin{equation*}
    = \frac{z}{(1-z)(1-yz)}. \quad \quad \quad \blacksquare
  \end{equation*}

\item [Study 2.] \textbf{2D first quadrant Upper VPV zeta function analogy for the previous exercise}. 

Show that for $\Re y > 1$, $\Re z > 1$,
    \begin{equation}\label{21.37a1}
\sum_{\substack{\gcd(a,b)=1 \\ 1 \leq a < b \leq 2}} \frac{1}{a^y b^z} = \frac{\zeta(y)\zeta(z)}{\zeta(y+z)}.
  \end{equation}
Proof: Starting with the left side of (\ref{21.37a1})
  \begin{eqnarray*}
    LHS  &=& \left( \frac{1}{1^y} \right) \frac{1}{2^z}\\
    &+& \left(\frac{1}{1^y} + \frac{1}{2^y} \right) \frac{1}{3^z} \\
    &+& \left(\frac{1}{1^y} + \frac{1}{3^y} \right) \frac{1}{4^z} \\
    &+& \left(\frac{1}{1^y}  + \frac{1}{2^y} + \frac{1}{3^y}  + \frac{1}{4^y} \right) \frac{1}{5^z}  \\
    &+& \left(\frac{1}{1^y}  + \frac{1}{5^y} \right) \frac{1}{6^z} \\
    &+& \left(\frac{1}{1^y}  + \frac{1}{2^y} + \frac{1}{3^y}  + \frac{1}{4^y} + \frac{1}{5^y} + \frac{1}{6^y} \right) \frac{1}{7^z}  \\
    &+& etc.
  \end{eqnarray*}
  \begin{equation*}
    = \sum_{\substack{\gcd(j,k)=1 \\ 0 \leq j < k \leq 1}} \sum_{a=0}^{\infty} \frac{1}{j^{ay} k^{az}}
    = \sum_{a=0}^{\infty} \sum_{\substack{\gcd(j,k)=1 \\ 0 \leq j < k \leq 1}}\frac{1}{j^{ay} k^{az}}
  \end{equation*}
    \begin{eqnarray*}
      &=& y^0 z^1 \\
      &+& (y^0 + y^1) z^2 \\
      &+& (y^0 + y^1 + y^2) z^3 \\
      &+& (y^0 + y^1 + y^2 + y^3) z^4  \\
      &+& (y^0 + y^1 + y^2 + y^3 + y^4) z^5 \\
      &+& (y^0 + y^1 + y^2 + y^3 + y^4 + y^5) z^6  \\
      &+& etc.
  \end{eqnarray*}
  \begin{equation*}
    = \frac{1-y}{1-y}z + \frac{1-y^2}{1-y}z^2 + \frac{1-y^3}{1-y}z^3 + \frac{1-y^4}{1-y}z^4 + \cdots
  \end{equation*}
  \begin{equation*}
    = \frac{1}{1-y} \left( \frac{z}{1-z} - \frac{yz}{1-yz} \right)
  \end{equation*}
  \begin{equation*}
    = \frac{z}{(1-z)(1-yz)}. \quad \quad \quad \blacksquare
  \end{equation*}
Infer from (\ref{21.37a}) and (\ref{21.37a1}) that they encode two similar statements about partitions into lattice point vectors $\langle a,b \rangle$ with $gcd(a,b)=1$ and for (\ref{21.37a}) we have non-negative integers $a$ and positive integers $b$; whereas for (\ref{21.37a1}) we have positive integers for both $a$ and $b$.

\item [Study 3.] \textbf{Find the particular cases of (\ref{21.37a})}. 

Prove the following:
\begin{eqnarray*}
\sum_{\substack{\gcd(a,b)=1 \\ 0 \leq a < b \leq 1}} \frac{1}{3^{n} \, 2^{n} -1} &=& \frac{2^{-n}}{(1-3^{-n})(1-2^{-n})} \quad (for \; \Re n >1), \\
   &=& \frac{1}{2^{n}}\left(1+\frac{1}{2^n}+\frac{1}{3^n}+\frac{1}{4^n}+\frac{1}{6^n}+\frac{1}{8^n}+\frac{1}{9^n}+\frac{1}{12^n}+\cdots  \right); \\
\sum_{\substack{\gcd(a,b)=1 \\ 0 \leq a < b \leq 1}} \frac{z^{a+b}}{1 - z^{a+b}} &=& \frac{z}{(1-z)^2}, \quad |z|<1; \\
\sum_{\substack{\gcd(a,b)=1 \\ 0 \leq a < b \leq 1}} tan^{2(a+b)}(\theta) &=& \frac{1}{cos^{2}(\theta)}.
\end{eqnarray*}

\item [Study 4.] \textbf{Find the particular cases of (\ref{21.37a1})}. 

Prove the following:
\begin{eqnarray*}
  \sum_{\substack{\gcd(a,b)=1 \\ a,b \geq 1}} \frac{1}{a^2 b^2} &=& \frac{5}{2}, \\
  \sum_{\substack{\gcd(a,b)=1 \\ a,b \geq 1}} \frac{1}{a^2 b^3} &=& \frac{\pi^2 \zeta(3)}{6 \zeta(5)}, \\
  \sum_{\substack{\gcd(a,b)=1 \\ a,b \geq 1}} \frac{1}{a^3 b^5} &=& \frac{9450\zeta(3)\zeta(5)}{\pi^8}.
\end{eqnarray*}

\item [Study 5.] \textbf{3D first hyperquadrant sums using gcd}. 

For each of $|x|, |y|, |z|<1,$ show that
  \begin{equation}\label{21.37b}
    \sum_{\substack{\gcd(a,b,c)=1 \\ a,b \geq 0,c>0}} \frac{x^a \, y^b \, z^c}{1 - x^a \, y^b \, z^c}
    = \frac{z}{(1-x)(1-y)(1-z)}.
  \end{equation}
    Similarly, show that for $\Re x > 1$, $\Re y > 1$, $\Re z > 1$,
    \begin{equation}\label{21.37b1}
    \sum_{\substack{\gcd(a,b,c)=1 \\ a,b,c \geq 1}} \frac{1}{a^x b^y c^z} = \frac{\zeta(x)\zeta(y)\zeta(z)}{\zeta(x+y+z)}.
  \end{equation}
Infer from (\ref{21.37b}) and (\ref{21.37b1}) that they encode two similar statements about partitions into lattice point vectors $\langle a,b,c \rangle$ with $gcd(a,b,c)=1$ and for (\ref{21.37b}) we have non-negative integers $a$, $b$ and positive integers $c$; whereas for (\ref{21.37b1}) we have positive integers for $a$, $b$ and $c$.

\item [Study 6.] \textbf{4D hyperquadrant VPV sums using gcd}. 

For each of $|w|, |x|, |y|, |z|<1,$ show that
  \begin{equation}\label{21.37c}
    \sum_{\substack{\gcd(a,b,c,d)=1 \\ a,b,c\geq0,d>0}} \frac{w^a \, x^b \, y^c \, z^d}{1 - w^a \, x^b \, y^c \, z^d}
    = \frac{z}{(1-w)(1-x)(1-y)(1-z)}.
  \end{equation}
      Similarly, show that for $\Re w$, $\Re x$, $\Re y$, and $\Re z$ all $ > 1$,
    \begin{equation}\label{21.37c1}
    \sum_{\substack{\gcd(a,b,c,d)=1 \\ a,b,c,d \geq 1}} \frac{1}{a^w b^x c^y d^z} = \frac{\zeta(w)\zeta(x)\zeta(y)\zeta(z)}{\zeta(w+x+y+z)}.
  \end{equation}
Infer from (\ref{21.37c}) and (\ref{21.37c1}) that they encode two similar statements about partitions into lattice point vectors $\langle a,b,c,d \rangle$ with $gcd(a,b,c,d)=1$ and for (\ref{21.37c}) we have non-negative integers $a$, $b$, $c$ and positive integers $d$; whereas for (\ref{21.37c1}) we have positive integers for $a$, $b$, $c$ and $d$.

\item [Study 7.] \textbf{5D hyperquadrant VPV sums using gcd}. 

For each of $|v|, |w|, |x|, |y|, |z|<1,$ show that
  \begin{equation}\label{21.37d}
    \sum_{\substack{\gcd(a,b,c,d,e)=1 \\ a,b,c,d \geq0,e>0}} \frac{v^a \, w^b \, x^c \, y^d \, z^e}{1 - v^a \, w^b \, x^c \, y^d \, z^e}
    = \frac{z}{(1-v)(1-w)(1-x)(1-y)(1-z)}.
  \end{equation}
      Similarly, show that for $\Re v$, $\Re w$, $\Re x$, $\Re y$, and $\Re z$ all $ > 1$,
    \begin{equation}\label{21.37d1}
    \sum_{\substack{\gcd(a,b,c,d,e)=1 \\ a,b,c,d,e \geq 1}} \frac{1}{a^v b^w c^x d^y e^z} = \frac{\zeta(v)\zeta(w)\zeta(x)\zeta(y)\zeta(z)}{\zeta(v+w+x+y+z)}.
  \end{equation}
Infer from (\ref{21.37d}) and (\ref{21.37d1}) that they encode two similar statements about partitions into lattice point vectors $\langle a,b,c,d \rangle$ with $gcd(a,b,c,d,e)=1$ and for (\ref{21.37d}) we have non-negative integers $a$, $b$, $c$, $d$ and positive integers $e$; whereas for (\ref{21.37d1}) we have positive integers for $a$, $b$, $c$, $d$ and $e$.

\item [Study 8.] \textbf{nD hyperquadrant sums using gcd}. 

For each of $|q_1|, |q_2|, |q_3|, \ldots, |q_n|<1,$ show that
  \begin{equation}\label{21.37e}
    \sum_{\substack{\gcd(a_1,a_2,a_3,\ldots,a_n)=1 \\ a_1,a_2,\ldots,a_{n-1} \geq 0,a_n \geq 1}} \frac{q_1^{a_1} \, q_2^{a_2} \, q_3^{a_3} \cdots  q_n^{a_n}}{1 - q_1^{a_1} \, q_2^{a_2} \, q_3^{a_3} \cdots q_n^{a_n}}
    = \frac{q_n}{(1-q_1)(1-q_2)(1-q_3)\cdots(1-q_n)}.
  \end{equation}
      Similarly, show that for $\Re q_1$, $\Re q_2$, $\Re q_3$, $\ldots$, and $\Re q_n$ all $ > 1$,
    \begin{equation}\label{21.37d1r}
  \sum_{\substack{\gcd(a_1,a_2,a_3, \ldots ,a_n)=1 \\ a_1,a_2 \ldots a_{n} \geq 1}} \frac{1}{a_1^{q_1} a_2^{q_2} a_3^{q_3} \cdots a_n^{q_n}} = \frac{\zeta(q_1)\zeta(q_2)\zeta(q_3) \ldots \zeta(q_n)}{\zeta({q_1}+ {q_2}+ {q_3}+ \cdots +{q_n})}.
  \end{equation}
Infer from (\ref{21.37d}) and (\ref{21.37d1}) that they encode two similar statements about partitions into lattice point vectors $\langle a_1,a_2,a_3,\ldots,a_n \rangle$ with $gcd(a_1,a_2,a_3,\ldots,a_n)=1$ and for (\ref{21.37d}) we have non-negative integers $a_1$, $a_2$, $\ldots$, $a_{n-1}$ and positive integers $a_n$; whereas for (\ref{21.37d1r}) we have positive integers for $a_1$, $a_2$, $\ldots$, up to $a_{n}$.

\end{itemize}

\end{document}